\documentclass[12pt]{article}

\usepackage[utf8]{inputenc}
\usepackage[T1]{fontenc}

\usepackage[francais]{babel}

\title{ENTROPIE SOFIQUE\\ d'apr\`es Lewis~Bowen, David Kerr et Hanfeng Li
\\
{\large S\'eminaire Bourbaki du 16 janvier 2016}}
\author{Damien GABORIAU}
\date{}

\usepackage{amsmath} 
\usepackage{amssymb}

\usepackage{color}

\usepackage{theorem}
\theoremstyle{changebreak}
\newtheorem{Theoreme}{Th\'eor\`eme}[section]

\newtheorem{prop}[Theoreme]{Proposition}

\newtheorem{theo}[Theoreme]{Th\'eor\`eme}
\newtheorem{defi}[Theoreme]{D\'efinition}
\newtheorem{lemm}[Theoreme]{Lemme}
\newtheorem{theodef}[Theoreme]{Th\'eor\`eme \&\ D\' efinition}
\newtheorem{coro}[Theoreme]{Corollaire}

 \usepackage[utf8]{inputenc}

\usepackage[bookmarks,bookmarksnumbered,pdfstartview=FitBH,backref=page,pdfpagelabels]{hyperref}
\hypersetup{
colorlinks=true, 
citecolor= blue,
linktoc= all, 
linkcolor= blue, 
urlcolor= green,
}

\begin{document}
\maketitle

\abstract{L'entropie, en systèmes dynamiques, fut introduite par A. Kolmogorov. Initialement focalisée sur les itérations d'une transformation préservant une mesure finie, la notion fut peu à peu généralisée, jusqu'à embrasser les actions des groupes moyennables ainsi que les actions topologiques. L. Bowen (2008) parvint à franchir la barrière du non moyennable en introduisant l'entropie sofique. Cet invariant rend les mêmes services que l'entropie classique pour les actions mesurées des groupes sofiques (une classe qui contient les groupes résiduellement finis). En 2010, D. Kerr et H. Li mirent au point une version topologique et un principe variationnel.
}

\medskip
Mots-clefs : systèmes dynamiques, théorie ergodique, entropie métrique, entropie topologique, entropie sofique, groupes sofiques

 


\vskip 30pt
\noindent{\bf INTRODUCTION}

L'entropie, dans les systèmes dynamiques, constitue une famille d'invariants multiformes avec des ramifications en grand nombre.
L'objectif de ce texte est de fournir une modeste introduction à un sujet qui a fleuri en bouquets depuis 2008 : l'entropie sofique à la suite des travaux de L.~Bowen, D.~Kerr et H.~Li.

On considérera l'entropie classique au sens de A.~Kolmogorov et Y.~Sina{\u\i}, de D.~Ornstein et B.~Weiss et l'entropie sofique au sens de L.~Bowen, D.~Kerr \& H.~Li  ainsi qu'une entropie dite de Rokhlin. Il s'agit dans tous ces cas d'attacher un nombre à une action d'un groupe dénombrable $G$, préservant une mesure $\mu$ de probabilité\footnote{On utilisera l'abréviation \textbf{p.m.p.} pour \og préservant la mesure de probabilité\fg.}, sur l'espace borélien standard $X$, c'est-à-dire à un homomorphisme de $G$ dans le groupe $\mathrm{Aut}(X,\mu)$ des bijections de $X$ préservant $\mu$, modulo égalité presque partout.

Ces notions seront ensuite mutées en leur variante topologique, c'est-à-dire pour des actions continues de $G$ sur un compact métrisable $X$, auxquelles elles seront reliées par un principe variationnel (travaux de D.~Kerr \& H.~Li).

Ces invariants fournissent des éléments de réponse au problème général suivant~: déterminer quand deux actions de $G$ sont conjuguées\footnote{Dans la catégorie mesurée ou dans la catégorie topologique.}.

Il s'agira dans cet exposé de décrire des résultats récents de Lewis Bowen, David Kerr et Hanfeng Li, qui étendent ces invariants  au-delà des groupes moyennables : jusqu'aux groupes dénombrables sofiques.

On n'y trouvera essentiellement pas de démonstrations, mais plutôt des éléments de comparaison et quelques références.
On a préféré parfois une présentation et des notations suggestives, évocatrices, et créant des échos entre différentes parties du texte, à d'autres plus strictes mais plus lourdes, en espérant rester dans des limites qui facilitent la lecture sans la rendre équivoque.

On a aussi fait le choix, certainement discutable, de repousser à la section~\ref{sect: definitions} un certain nombre de définitions standard (partition génératrice, décalage,...) et à la section~\ref{sect:groupes sofiques} des explications concernant les groupes sofiques,
afin de permettre de rentrer plus directement dans le vif du sujet. Le lecteur peu familier avec ces notions pourra s'y reporter.

{Ce texte accompagne mon exposé au séminaire Bourbaki du 16 janvier 2016, accessible en ligne~:
\url{https://www.youtube.com/watch?v=TW8O9Lo631w}.}

\section{Présentation de la thématique}

L'entropie, en théorie de l'information, est un concept fondamental introduit par C.~Shannon en 1948 \cite{Shannon-1948-entropy}. 
A.~Kolmogorov \cite{Kolmogorov-1958-entropy,Kolmogorov-1959-entropy-erratum} l'a utilisée pour définir un invariant non moins fondamental en théorie ergodique~: l'entropie d'une transformation préservant la mesure, dont il a posé les bases entre 1958 et 1962 avec quelques proches mathématiciens, notamment Y.~Sina{\u\i} et V.~Rokhlin.

L'entropie de Shannon d'une partition dénombrable $\alpha=(A_i)_{i\in K}$ de l'espace de probabilité $(X,\mu)$ est définie par 
\begin{equation}
H(\alpha)\overset{\mathrm{def}}{=}-\sum_{i\in K} \mu(A_i) \log(\mu(A_i)). \tag{\textbf{Entropie de Shannon}}
\end{equation}

\subsection{Pour le groupe \texorpdfstring{$G=\mathbf{Z}$}{Z}}
Pour un isomorphisme préservant la mesure  $T\colon X\to X$, ou autrement dit une action \textbf{p.m.p.} $G\curvearrowright^{\!T}{\!}(X,\mu)$ de $G = \mathbf{Z}$, sur l'espace de probabilité standard, la définition de Kolmogorov nécessite l'existence d'une partition \textbf{génératrice}\footnote{Avec l'aide de l'action de $G$, elle permet de séparer presque tous les points de $X$ ; voir 
section~\ref{sect: definitions}.} $\alpha$ d'entropie de Shannon finie et considère les entropies de Shannon des joints des itérés (les raffinements de $\alpha$ obtenus par itérations) et normalisées~:
\begin{equation}
h({G} \curvearrowright^{\!T}{\!} X,\mu)\overset{\mathrm{def}}{=}\lim\limits_{n\to\infty} \frac{1}{n} H(\vee_{i=0}^{n} T^{-i}\alpha).\tag{\textbf{Entropie de Kolmogorov}}
\end{equation}
Le point clef est bien entendu l'indépendance vis-à-vis de la partition génératrice (l'existence de la limite n'est pas difficile).

Y.~Sina{\u\i} \cite{Sinai-1959-entropy} a apporté une amélioration significative en observant que, parmi les partitions d'entropie de Shannon finie, les partitions génératrices, lorsqu'elles existent, maximisent la quantité considérée (c'est ce qu'on appelle le {\em théorème de Kolmogorov-Sina{\u\i}}).
Cela permet de se débarrasser de l'hypothèse génératrice~:
\begin{equation}
h_{\mathrm{KS}}({G} \curvearrowright^{\!T}{\!} X,\mu)\overset{\mathrm{def}}{=}\sup_{\substack{\beta \textrm{ partition } \\ H(\beta)<\infty}} \lim\limits_{n\to\infty} \frac{1}{n} H(\vee_{i=0}^{n} T^{-i}\beta).
\tag{\textbf{Entropie de Kolmogorov-Sina{\u\i}}}
\end{equation}

L'un des tout premiers succès de cette théorie aura été de résoudre le problème, qu'on attribue à J.~von~Neumann, de la conjugaison mesurée des actions par décalage de Bernoulli sur des espaces de base à $2$ points $G\curvearrowright \{1,2\}^G$, respectivement $3$ points $G\curvearrowright \{1,2,3\}^G$ (avec les mesures canoniques). Leur entropie valant $\log(2)$, resp. $\log(3)$, ils ne peuvent pas être isomorphes.
Et on touche  de suite à deux propriétés cruciales de l'entropie de Kolmogorov~:
\begin{itemize}
\item [\textbf{(a)}]  l'entropie est un \textbf{invariant de conjugaison mesurée}, 
\item [\textbf{(b)}] l'entropie des \textbf{décalages de Bernoulli} est  égale à l'entropie de Shannon de leur base~: 
$h_{\mathrm{KS}}({G}\curvearrowright (L,\lambda)^{{G}})= H(\lambda) \overset{\mathrm{def}}{=} \sum_{l\in L} -\lambda(l) \log\lambda(l).$

En particulier, une conjugaison mesurée entraîne que les bases ont même entropie~:
\begin{equation}
G\curvearrowright(K,\kappa)^{G}\simeq G\curvearrowright(L,\lambda)^{G} \Rightarrow  H(\kappa)=H(\lambda). \tag{\textbf{Théorème de Kolmogorov}}\end{equation}
\end{itemize}
Ce résultat reste vrai, que la base (et son entropie de Shannon) soit finie ou non\footnote{L'entropie de Shannon d'un espace diffus tel que $[0,1]$ muni la mesure Lebesgue vaut $+\infty$.}. Voir aussi section~\ref{sect: ent Rokhlin} pour l'existence de partitions génératrices finies.

Après Kolmogorov, le problème est devenu celui d'une réciproque. D.~Ornstein a développé de puissantes méthodes qui lui ont permis d'identifier une forme d'ubiquité des actions Bernoulli (on peut consulter l'article de survol \cite{Ornstein-2013-survey-Bernoulli}). Et la réciproque en a découlé \cite{Ornstein-1970-entrop-shift-isom}~:
\begin{equation}
H(\kappa)=H(\lambda) \Rightarrow G\curvearrowright(K,\kappa)^{G}\simeq G\curvearrowright(L,\lambda)^{G}. \tag{\textbf{Théorème d'Ornstein}}
\end{equation}
Ainsi, l'entropie est un invariant complet de conjugaison mesurée parmi les décalages de Bernoulli.

Une propriété saute aux yeux avec le théorème de Kolmogorov-Sina{\u\i} (elle lui est d'ailleurs essentiellement équivalente), c'est le comportement de l'entropie sous \textbf{facteur}\footnote{Les facteurs jouent un grand rôle en théorie ergodique ; ils sont liés aux partitions non génératrices (voir section~\ref{sect: definitions}).}, \emph{i.e.} une application mesurable $\pi\colon (X,\mu)\to (Y,\nu)$ qui pousse $\mu$ sur $\nu$ (essentiellement surjective, donc -- l'image est de mesure pleine) et qui est équivariante pour des actions  ${G}\curvearrowright(X,\mu)$ et ${G}\curvearrowright(Y,\nu)$. C'est une troisième propriété fondamentale de l'entropie :
\begin{itemize}
\item[\textbf{(c)}] l'entropie décroît sous \textbf{facteur}~: $h_{\mathrm{KS}}({G}\curvearrowright X,\mu)\geqslant h_{\mathrm{KS}}({G}\curvearrowright Y,\nu)$.
\end{itemize}
En effet, les partitions de $Y$ et les calculs d'entropie de joints se remontent dans $X$.

Signalons une forme de réciproque due à Y.~Sina{\u\i}.
\\
\emph{Si $\mathbf{Z} \curvearrowright^{\! T}\! (X,\mu)$ est une action p.m.p. ergodique, alors elle factorise sur tout décalage de Bernoulli $\mathbf{Z}\curvearrowright (K^{\mathbf{Z}}, \kappa^{\mathbf{Z}})$ à base dénombrable \og entropie compatible \fg, \emph{i.e.} qui vérifie $h_{\mathrm{KS}}(\mathbf{Z} \curvearrowright^{\! T}\! X,\mu)\geqslant h_{\mathrm{KS}}(\mathbf{Z} \curvearrowright K^{\mathbf{Z}}, \kappa^{\mathbf{Z}})$}  \cite{Sinai-1962-weak-isom}.
\label{th: Sinai fact sur Bernoulli}

Les applications de l'entropie d'une transformation se sont diffusées dans une grande partie des systèmes dynamiques, et les résultats sont quasiment indénom\-brables.
Nous reviendrons plus loin sur certains d'entre eux. Pour un survol historique détaillé de l'entropie de Kolmogorov, on pourra consulter l'article \cite{Katok-2007-50-yrs-entropy}.

\subsection{Pour les groupes $G$ moyennables}
La théorie s'est également développée pour des actions de groupes plus généraux que le groupe $\mathbf{Z}$. Et on s'est vite aperçu que les choses se passaient bien pour les groupes commutatifs, puis pour le dire rapidement, la théorie s'est étendue à tous les groupes dénombrables moyennables (voir section~\ref{sect:groupes sofiques}), avec le considérable travail fondateur de D.~Ornstein et B.~Weiss \cite{OW87}.
Une particularité de la définition de l'entropie de Kolmogorov pour $G=\mathbf{Z}$ est l'utilisation des intervalles $\{0, 1, \cdots, n\}$ qui ont un petit bord (formé de deux points $0$ et $n$).
La définition de l'entropie pour les actions p.m.p. des groupes moyennables doit maintenant faire appel aux \textbf{suites de F{\o}lner}\footnote{Des parties finies à constante isopérimétrique tendant vers $0$, caractéristiques de la moyennabilité de $G$.} 
 dans le groupe $G$, mais l'essentiel de la théorie est conservé, et principalement les propriétés (a) (conjugaison mesurée), (b) (décalages de Bernoulli) et (c) (décroissance sous facteur) mises en évidence ci-dessus. Le théorème d'Ornstein reste également valide et l'entropie mesurée est un invariant complet de conjugaison mesurée parmi les décalages de Bernoulli.

\subsection{Pour les groupes $G$ non moyennables}
L'extension de la théorie entropique à des groupes non moyennables est restée une quête inaccessible jusqu'à l'intervention en 2008 de L.~Bowen et de son entropie sofique.
Une quête d'autant plus frustrante que A.~Stepin \cite{Stepin-1975} a montré que le théorème d'Ornstein, la partie réciproque donc, s'étendait facilement des sous-groupes au groupe ambiant. Ainsi, 
pour tous les groupes contenant un groupe moyennable infini\footnote{On peut signaler que L.~Bowen \cite{Bowen-2012-almost-Ornstein}  a étendu cette propriété 
\og si $H(\nu_1)=H(\nu_2)$, alors  $G\curvearrowright (K_1^G,\nu_1^G)\simeq G\curvearrowright (K_2^G,\nu_2^G)$\fg\ à tous les groupes infinis dénombrables, à condition qu'aucun des deux espaces de base ne soit constitué de seulement deux atomes.}~: 
{\em si les entropies de Shannon des espaces de base sont égales, alors les décalages de Bernoulli associés sont conjugués}, sans qu'on ait été capable de dire si
réciproquement, comme dans le théorème de Kolmogorov, la conjugaison entraînait l'égalité des entropies de Shannon des bases.

\subsection{Le blocage psychologique de l'exemple d'Ornstein-Weiss}
\label{sect: appl Ornstein-Weiss}
Dans leur article \cite{OW87}, D.~Ornstein et B.~Weiss ont exhibé un exemple qui semblait rendre vaine toute tentative au-delà du moyennable\footnote{Rappelons que les prototypes de groupes non moyennables sont les groupes contenant le groupe libre $\mathbf{L}_2$.}.
Soit  $\mathbf{L}_2=\langle a, b\rangle$ le groupe libre à deux générateurs et soit $\mathbb{K}$ un groupe fini abélien, par exemple $\mathbb{K}=\mathbf{Z}/2\mathbf{Z}$.
L'application 
\begin{equation}
\label{eq: appl OW}
\Theta \colon \left(\begin{array}{cclcc}
 \mathbb{K}^{\mathbf{L}_2}&\to &\mathbb{K}^{\mathbf{L}_2}\times \mathbb{K}^{\mathbf{L}_2}\simeq (\mathbb{K}\times \mathbb{K})^{\mathbf{L}_2}\\
\omega &\mapsto &
(\omega(g a)-\omega(g), \omega(g b)-\omega(g) )_{g\in \mathbf{L}_2}
\end{array}\right)
\end{equation}
est surjective et définit une factorisation du décalage de Bernoulli de base $\mathbb{K}$ sur celui de base $\mathbb{K}\times \mathbb{K}$, où les bases sont munies des mesures d'équiprobabilité.
Toute bonne théorie devant leur donner entropie $\log \vert \mathbb{K}\vert $, resp. $\log \vert \mathbb{K}\times \mathbb{K} \vert=2 \log \vert \mathbb{K}\vert  $, on aurait affaire à {\em un facteur qui augmente l'entropie} !

Et cette application est loin d'être une pathologie, c'est un homomorphisme continu $\mathbf{L}_2$-équivariant de groupes compacts, qu'on équipe de leurs mesures de Haar. Son noyau est fini~: ce sont les fonctions constantes $\mathbf{L}_2\to \mathbb{K}$. 
On peut aussi l'interpréter comme l'application cobord qui va des cochaînes de dimension $0$ à coefficients dans $\mathbb{K}$ dans celles de dimension $1$, pour l'arbre de Cayley $\mathcal{T}$ du groupe libre.
Quelques détails supplémentaires sont fournis en section~\ref{sect:plus sur l'ex de OW}.

 Cet exemple a conduit D.~Ornstein et B.~Weiss \cite{OW87} à demander si tous les décalages de Bernoulli sur un groupe non moyennable étaient isomorphes\footnote{Auquel cas, l'entropie de Shannon de la base n'aurait pas été un invariant de conjugaison.}.
Et pour enfoncer le clou\footnote{Le théorème \og entropie compatible \fg\ de Sina{\u\i} (section \ref{th: Sinai fact sur Bernoulli}) qui montrait l'isomorphisme {\em faible} des {$\mathbf{Z}$-décalages} de même entropie, fut considéré comme un jalon important vers le théorème d'Ornstein d'isomorphisme {\em fort} de ces décalages.}, 
L. Bowen montre que~:
\begin{theo}[Bowen {\cite[Th. 1.3]{Bowen-2011=weak-isom-Bernoulli}}]\label{th:Bowen isom faibles des Bernoulli}
Pour tout groupe $G$ contenant $\mathbf{L}_2$, tous les décalages de Bernoulli factorisent les uns sur les autres (on dit qu'ils sont faiblement isomorphes).
\end{theo}
C'est une catastrophe ?!
\`A moins de renoncer à la propriété de décroissance... et de changer de point de vue.

\subsection{Point de vue externe, modèles finis}
L.~Bowen dans une série de résultats retentissants annoncés à partir de 2008 va faire rebondir complètement le sujet.
Le premier de la série \cite{Bowen-2010-f-invariant} (annoncé en 2008, publié en 2010) traite spécifiquement le cas où $G$ est un groupe libre. Il introduit une quantité qu'il appelle le $f$-invariant\footnote{Le symbole $f$ dans $f$-invariant est utilisé pour évoquer le qualificatif  {\em free} de {\em free group}.} dont il montre qu'il est invariant de conjugaison mesurée et qu'il prend la valeur $H(\kappa)$ sur un décalage de Bernoulli $G\curvearrowright (K,\kappa)^G$ dès que $K$ est fini\footnote{Il parvient aussi à traiter des $K$ infinis d'entropie de Shannon finie.}.

Dans cet article et dans le suivant \cite{Bowen-2010-mes-conj-inv-sofic}, il adopte un \og point de vue externe\fg, comme dit D.~Kerr. Il cherche à modéliser sur des ensembles finis $D$ la dynamique de l'action $G\curvearrowright^{\!T}{\!} (X,\mu)$. 
Le point de vue adopté est plus proche de celui de l'interprétation statistique de l'entropie de Boltzmann. L'entropie de Shannon $H(\alpha)$ d'une partition $\alpha=(A_k)_{k\in K}$ de $X$ (sans considérer de dynamique)
peut s'obtenir de la façon suivante. 
On considère un ensemble fini $D$ (qu'on pense très grand) et on regarde toutes les partitions sur $D$ qui imitent bien $\alpha$ (au sens où les pièces ont les bonnes mesures à $\epsilon$ près) ; on estime leur nombre, puis on prend le taux de croissance exponentielle de ce nombre en la taille de $D$ qu'on fait tendre vers l'infini, puis on serre les $\epsilon$. Plus précisément, c'est une forme simple d'un principe de grande déviation~:

\begin{prop}[Boltzmann, Sanov {\cite{Sanov-1957-large-dev}}]
\label{prop:Boltzmann-Sanov}
\begin{equation*}
H(\alpha)=\lim\limits_{\epsilon\to 0} \lim\limits_{\vert D\vert \to \infty} \frac{1}{\vert D\vert }\log \left\vert 
\left\{\textrm{partitions } (V_k)_{k\in K} \textrm{ de } D \ \Big\vert \ \forall k\in K,\  \left\vert \frac{\vert V_k \vert}{\vert D\vert} -\mu(A_k)\right\vert <\epsilon\right\}
\right\vert.
\end{equation*}
\end{prop}
C'est ce point de vue, étendu à la situation où un groupe agit, qui conduit à la notion d'entropie sofique.

On revient en détail sur cet énoncé (en section~\ref{sect: preuve du lem Boltzmann-Sanov})
et on en donne une preuve \og probabiliste\fg\ qui reviendra en écho lors des estimations pour les décalages de Bernoulli (section~\ref{sect: ent sofique Bernoulli}).

\subsection{Introduction de la dynamique, soficité} 
Si on cherche maintenant à modéliser la dynamique de $G\curvearrowright^{\!T}{\!} (X,\mu)$ sur des ensembles finis $D$ qui auront vocation à devenir grands, il faudra commencer par concéder que le groupe $G$ possède des actions sur de tels ensembles, c'est-à-dire des homomorphismes \og pas trop triviaux\fg\ dans des groupes symétriques $\mathrm{Sym}(D)$. Un bon cadre est celui des groupes résiduellement finis (voir section~\ref{sect:groupes sofiques}). Mais puisque les erreurs sous-exponentielles seront tolérables\footnote{Voire, en vérité, des erreurs exponentielles mais de taux de croissance arbitrairement petit.}, on sera en mesure d'admettre de petites erreurs dans les homomorphismes, c'est-à-dire considérer des groupes sofiques.

Les \textbf{groupes sofiques} (introduits par M.~Gromov) sont des groupes  qui possèdent de bonnes propriétés d'approximation en termes de permutations sur des ensembles finis ; et un peu plus précisément, des suites de {\em presque-homomorphismes-sans-point-fixe} $\Sigma=(\sigma_n \colon G\to \mathrm{Sym}(D_n))_{n}$  dans des groupes de permutations (voir section~\ref{sect:groupes sofiques}).

\begin{defi}[Groupe sofique, approximation sofique] \label{defn: gp approx sofic}
Un groupe dénombrable $G$ est \textbf{sofique} s'il admet une suite 
d'ensembles finis $D_n$ et d'applications\footnote{Insistons : pas nécessairement des homomorphismes !}
${\sigma_n \colon G \rightarrow \mathrm{Sym}(D_n)}$  telles que $\sigma_n(1_G) = 1_{\mathrm{Sym}(D_n)}$ et
\begin{enumerate}
\item[\rm (i)] (\textbf{presque action}) \hfil
$\frac{1}{\vert D_n\vert} \vert \{ z \in D_n \ \vert \ \sigma_n(g) \circ \sigma_n(h)(z) = \sigma_n(g h)(z)\} \vert \underset{n\to \infty}{\longrightarrow} 1$, \hfil $\forall g, h \in G$
\item[\rm (ii)] (\textbf{presque libre}) \hfil
$\frac{1}{\vert D_n\vert} \vert\{ z \in D_n \ \vert \ \sigma_n(g)(z) \neq z\}\vert \underset{n\to \infty}{\longrightarrow} 1$, \hfil $\forall g\in G\setminus\{1_G\}$.
\end{enumerate}
Une suite $\Sigma=(\sigma_n \colon  G \rightarrow \mathrm{Sym}(D_n))_{n}$ comme ci-dessus, où le cardinal $d_n\overset{\mathrm{def}}{=}\vert D_n\vert$ tend vers l'infini\footnote{Cette condition, automatiquement satisfaite si le groupe $G$ est infini, permet d'éviter certaines pathologies pour les groupes finis.}   
est appelée une \textbf{approximation sofique} du groupe $G$. 
\end{defi}
Les groupes résiduellement finis rentrent dans cette classe grâce à leurs actions par multiplication sur leurs quotients finis. On peut déjà goûter tout le sel de la notion d'entropie sofique en se restreignant aux actions de ces derniers et on peut lire ce texte en se concentrant sur eux.

L'extension aux groupes sofiques n'est pas que pur désir de généralisation maximale. On aurait deux types de théories de l'entropie avec une intersection non triviale, l'une pour les groupes moyennables, l'autre pour les groupes résiduellement finis. Or, les groupes moyennables sont également sofiques et, pour eux, les invariants obtenus dans le cadre de cette théorie de l'entropie sofique, quoique définis de manière bien différente, se trouvent coïncider avec ceux de la théorie classique de l'entropie \cite{Bowen-2012-sofic-ent-amenab, Kerr-Li-2013-sofic-amenabl-dyn-entrop}.
On disposera ainsi d'un emboîtement de généralisations successives de l'entropie mesurée.

On développe un peu cette notion de soficité et on la relie à la moyennabilité et aux groupes résiduellement finis dans la section~\ref{sect:groupes sofiques}. On signale (voir Théorème~\ref{th: Kerr-Li entrop Gottschalk conject}) la preuve entropique due à D.~Kerr et H.~Li de la conjecture de surjonctivité de Gottschalk pour les groupes sofiques.

\subsection{Entropie sofique mesurée}
Considérons donc un groupe sofique $G$ et une de ses approximations sofiques  $\Sigma=(\sigma_n \colon G\to \mathrm{Sym}(D_n))_{n}$ (définition~\ref{defn: gp approx sofic}).
Les ensembles finis $D_n$ de $\Sigma$,  équipés de leur {\em presque action}\footnote{Comme dit précédemment, on peut faire semblant de croire que les $\sigma_n$ sont de vrais homomorphismes.} de $G$, sont munis chacun de leur mesure de probabilité uniforme $\mathbf{u}_{n}$. Ils sont envisagés comme des espaces modèles pour les diverses dynamiques produites par les actions de $G$. 

Considérons une action p.m.p. $G\curvearrowright^{\!T}{\!} (X,\mu)$. 
Soit $\alpha \colon X\to K$ une $K$-partition mesurée finie de $X$ (les pièces sont les $A_k=\alpha^{-1}(k)$, indexées, ou colorées si on veut, par l'ensemble fini $K$). L'ensemble $K^{D_n}=\{a \colon D_n\to K\}$ est alors simplement l'ensemble de toutes les $K$-partitions de $D_n$.

Si $F\subset G$ est une partie finie contenant $\mathrm{id}_G$, alors par itérations, elle définit la partition $F$-raffinée $\alpha^{\vee F}$ et de manière analogue, chaque partition $a \colon D_n\to K$ fournit une partition $F$-raffinée\footnote{Observons que la définition de $a^{\vee F}$ ne nécessite pas que $\sigma_n$ soit un homomorphisme. Une approximation sofique fera aussi bien l'affaire.} $a^{\vee F}$, 
\begin{equation}
\alpha^{\vee F} \colon  \left(\begin{array}{ccl}
X &\to &K^F 
\\
x &\mapsto & (\alpha(T(f)\cdot x))_{f\in F}
\end{array}\right)
\ \textrm{ et } \ 
a^{\vee F} \colon  \left(\begin{array}{ccl}
D_n&\to &K^F 
\\
{v} &\mapsto & (a(\sigma_n(f)\cdot {v}))_{f\in F}
\end{array}\right).
\end{equation}
Il s'agit alors de comparer les mesures des pièces de ces partitions dans $X$ et dans $D_n$.
Ou autrement dit de comparer les mesures poussées en avant $\alpha^{\vee F}_{*} \mu$ et $a^{\vee F}_{*} \mathbf{u}_n$ sur $K^F$.
On note 
\begin{equation}\label{eq:def M mu-intro}
\mathcal{M}_\mu(\alpha, F, \epsilon, \sigma_n)\overset{\mathrm{def}}{=}\left\{a \in K^{D_n}\ \big\vert \ \Vert \alpha^{\vee F}_{*} \mu-a^{\vee F}_{*} \mathbf{u}_n\Vert_{1} < \epsilon\right\}.
\end{equation}
C'est l'ensemble de toutes les $K$-partitions de $D_n$ qui, une fois itérées par $F\subset G$, fournissent des pièces de mesure proche de celle des pièces obtenues de la même façon pour $\alpha$ sur $X$, et la proximité est mesurée dans la norme $\ell^1$. Ceux qui préfèrent la version ensembliste pourront se reporter à la section~\ref{sect: modeles finis de la dynamique}.

L.~Bowen considère alors le cardinal de cet ensemble et son taux de croissance exponentielle en le cardinal de $D_n$, puis il minimise selon les paramètres $F$ et $\epsilon$. 
Le résultat frappant est qu'on obtient une quantité qui est indépendante du choix d'une partition génératrice finie (si elle existe !). 

\begin{theodef}[Entropie sofique mesurée, Bowen {\cite{Bowen-2010-mes-conj-inv-sofic}}]
\label{th-def: ent sof mes Bowen}
Soit $\Sigma$ une approximation sofique du groupe sofique $G$. Soit $G\curvearrowright^{\!T}{\!} (X,\mu)$ une action p.m.p. qui admet une partition génératrice finie $\alpha$. Alors, la quantité ci-dessous est indépendante du choix de la partition génératrice finie.
Cette valeur commune est appelée l'\textbf{entropie sofique mesurée} de l'action relativement à $\Sigma$, et on la note~:
\begin{equation}
h_{\mathrm{mes}}^{\Sigma}(G\curvearrowright^{\!T}{\!} X,\mu)\overset{\mathrm{def}}{=}\inf_{\epsilon > 0} \inf_{\substack{F \subset G\\F \text{ fini}}} \limsup_{n \rightarrow \infty} \frac{\log \left| \mathcal{M}_\mu(\alpha, F, \epsilon, \sigma_n) \right |}{\vert D_n \vert}.
\tag{\textbf{Entropie sofique mesurée}}\label{ent sofic Bowen}
\end{equation}
\end{theodef}
L'entropie sofique prend ses valeurs dans $[0, +\infty]\cup \{-\infty\}$. La valeur $-\infty$ correspond aux situations où l'ensemble $\mathcal{M}_\mu(\alpha, F, \epsilon, \sigma_n)$ est vide\footnote{Via la convention $\log \vert \emptyset\vert=-\infty$.} à partir d'un certain rang $n$ pour $F$ grand et $\epsilon$ petit. La valeur $+\infty$ apparaîtra lorsqu'on se sera débarrassé de l'hypothèse génératrice.
On peut en effet noter que si $\alpha  \colon X\to K$ est une partition génératrice finie, alors $\mathcal{M}_\mu(\alpha, F, \epsilon, \sigma_n)\subset K^{D_n}$ nous donne  $h_{\mathrm{mes}}^{\Sigma}(G\curvearrowright X,\mu) \leqslant  \frac {\log\vert K^{D_n} \vert }{\vert D_n \vert}=  \log\vert K \vert$ ; c'est un premier pas vers l'entropie de Rokhlin (section~\ref{sect: ent Rokhlin}).

La dépendance en l'approximation sofique $\Sigma$ est l'objet de grandes spéculations.
On dispose d'exemples où elle en dépend. Voir la section~\ref{sect: depend-approx-sofic} pour des informations plus détaillées.
Mais en l'état actuel des connaissances, tous ces exemples reposent de façon plus ou moins directe sur l'utilisation de la propriété $(\tau)$, une variante de la propriété (T) de Kazhdan, et sur le fait que pour certaines actions de certains groupes, on peut trouver des approximations sofiques $\Sigma_1$ qui donnent lieu à des  $\mathcal{M}_\mu(\alpha, F, \epsilon, \sigma_n)$ qui auront tendance à être vides (et donneront $h_{\mathrm{mes}}^{\Sigma_1}(G\curvearrowright X,\mu)=-\infty$), tandis que pour d'autres approximations $\Sigma_2$, ils ne seront pas vides.
On ne dispose pas d'exemple d'action pour laquelle deux approximations sofiques donneraient des valeurs réelles distinctes.
Il découle de cette dépendance que la $ \limsup_{n \rightarrow \infty}$ qui intervient dans la définition de l'entropie sofique ne peut pas être remplacée par une limite\footnote{
Penser à une approximation sofique qui \emph{piocherait} alternativement dans $\Sigma_1$ ou $\Sigma_2$.} et ce détail complique très sérieusement la vie.
L.~Bowen \cite[Rem. 1] {Bowen-2010-mes-conj-inv-sofic} signale que sa théorie peut aussi se développer de manière satisfaisante en remplaçant cette $\limsup$ par une $\liminf$ voire par une ultra-limite\footnote{Een introduisant un ultra-filtre comme paramètre supplémentaire. Et ce point de vue serait assez cohérent avec l'approche des groupes sofiques par les ultra-produits développée par G.~Elek et E.~Szab{\'o} dans \cite{2005=Elek-Szabo=hyperlinearity}.}. Une variante sans conséquence fondamentale consisterait à remplacer la norme $\ell^1$ par d'autres normes pour les mesures sur $K^{F}$.

L'entropie sofique vérifie, elle aussi, les propriétés (a) et (b).
\\ (a) C'est un invariant de conjugaison mesurée.
\\ (b) L'entropie sofique des décalages est bien celle de leur base~:
\begin{theo}[Bowen {\cite{Bowen-2010-mes-conj-inv-sofic}}]
 Si $G\curvearrowright (K^{G}, \nu^{\otimes G})=(K,\nu)^{G}$ est une action par décalage de Bernoulli d'un groupe sofique, alors pour toute approximation sofique $\Sigma$ de $G$, on a~:
\[h_{\mathrm{mes}}^{\Sigma}(G\curvearrowright K^{G}, \nu^{\otimes G}) =H(\nu).\tag{\textbf{Entropie des Bernoulli}}\]
\end{theo}
\`A la différence de l'entropie de Kolmogorov pour les groupes moyennables, la dé\-mons\-tration de cet énoncé n'est pas immédiate. Il est dû à L.~Bowen \cite{Bowen-2010-mes-conj-inv-sofic} lorsque l'entropie de la base est finie et à D.~Kerr et H.~Li~\cite{Kerr-Li-2011-Bernoulli-infinite-entropy} en général (en anticipant un peu sur la suite de ce texte). On trouvera des indications de preuve en section~\ref{sect: ent sofique Bernoulli}.
En combinant cela avec les résultats d'Ornstein, Stepin et Ornstein-Weiss on obtient~:
\begin{theo}
Pour tout groupe sofique $G$ contenant un sous-groupe moyennable infini\footnote{Ou bien 
$G$ sofique et aucun des espaces de base n'est constitué de seulement deux atomes \cite{Bowen-2012-almost-Ornstein}.}, l'entropie de Shannon de la base est un invariant complet de conjugaison mesurée~:
\begin{equation}G\curvearrowright(K,\kappa)^{G}\simeq G\curvearrowright(L,\lambda)^{G} \Longleftrightarrow  H(\kappa)=H(\lambda). \tag{\textbf{Invariant complet}}\end{equation}
\end{theo}
On ignore à ce jour si c'est vrai pour tout groupe dénombrable.

En revanche, la propriété (c) est maintenant mise en défaut notamment par l'application d'Ornstein-Weiss (section~\ref{sect: appl Ornstein-Weiss}).
{\em L'entropie sofique peut croître sous facteurs.}

\subsection{Entropie sofique sans partition génératrice}
\label{subsect: Ent sof sans part gen}
Et l'histoire semble se répéter.
Nous voici avec une notion d'entropie bien définie lorsqu'on dispose d'une partition génératrice finie. 
On aimerait bien se débarrasser de cette hypothèse, notamment dans la perspective d'un principe variationnel (voir section~\ref{sect: ent topologique Kerr-Li}).
Ce souhait sera réalisé par les travaux de D.~Kerr et H.~Li  \cite{Kerr-Li-2011-Variationnal-principle, Kerr-Li-2013-sofic-amenabl-dyn-entrop}. Signalons que L.~Bowen, dans son article fondateur  \cite{Bowen-2010-mes-conj-inv-sofic}, étend par un procédé limite une partie de ses résultats aux actions qui possèdent une partition génératrice dénombrable d'entropie de Shannon finie.

D.~Kerr et H.~Li développent une théorie de l'entropie sofique mesurée dans un cadre d'analyse fonctionnelle \cite{Kerr-Li-2011-Variationnal-principle} (où l'hypothèse de génération finie perdra de sa pertinence) et démontrent qu'elle est équivalente à celle de L.~Bowen en présence d'une partition génératrice finie.
 Ils introduisent une nouvelle notion, celle d'\textbf{entropie sofique topologique} et démontrent un \textbf{principe variationnel} (voir section~\ref{sect: ent topologique Kerr-Li}).

Ils parviennent ensuite à réintégrer ces notions dans un cadre analogue à celui de L.~Bowen de partitions finies et de dénombrements dans \cite{Kerr-2013-Sofic-meas-ent-via-finite-partitions} pour l'entropie mesurée et dans \cite{Kerr-Li-2013-ent-top-version-combinat} pour l'entropie topologique. C'est plutôt sur ces dernières versions qu'on va se concentrer.

Partitions non génératrices et facteurs étant intimement liés (voir section~\ref{subsubsect:partitions generatrices}), l'exemple d'Ornstein-Weiss, ainsi que le théorème~\ref{th:Bowen isom faibles des Bernoulli} de factorisation les uns sur les autres des décalages de Bernoulli des groupes contenant $\mathbf{L}_2$ \cite[Th. 1.3]{Bowen-2011=weak-isom-Bernoulli}, indiquent qu'une définition à la Sina{\u\i}, en prenant un \emph{supremum} sur toutes les partitions finies de l'entropie sofique introduite par Bowen, conduirait certainement à un invariant peu intéressant (qui donnerait par exemple la valeur $+\infty$ pour tous les décalages des groupes contenant $\mathbf{L}_2$). 

L'idée est alors de définir l'entropie sofique d'une partition finie $\beta  \colon  X \rightarrow L$, à valeurs dans l'ensemble fini $L$, en la confrontant à toutes les partitions finies mesurables $\alpha  \colon  X \rightarrow K$ qui sont \textbf{plus fines}\footnote{Les pièces de $\beta$ sont obtenues en regroupant entre elles des pièces de $\alpha$.} qu'elle (ce qu'on note  $\alpha \geqslant \beta$), c'est-à-dire telles que $\beta = \Theta_{\beta,\alpha} \circ \alpha$ pour une certaine application  $\Theta_{\beta,\alpha}  \colon  K \rightarrow L$ (de fusion des pièces).
\begin{equation}
\left.\begin{array}{rcl} X & \overset{\alpha}{\longrightarrow}  & K \\ & {\searrow}{\beta} & \downarrow \Theta_{\beta,\alpha} \\ &  & L\end{array}\right.
\left.\begin{array}{rcl} D_n & \overset{a}{\longrightarrow}  & K \\ & {\searrow}{b} & \downarrow \Theta_{\beta,\alpha} \\ &  & L\end{array}\right.
\tag{\textbf{Fusion des pièces}}
\end{equation}
Les partitions de $D_n$ qui sont de bons $K$-modèles finis de $\alpha$ pour $(F,\epsilon)$, c'est-à-dire les éléments de 
$\mathcal{M}_\mu(\alpha, F, \epsilon, \sigma_n)=\left\{a \in K^{D_n}\ \big\vert \ \Vert \alpha^{\vee F}_{*} \mu-a^{\vee F}_{*} \mathbf{u}_n\Vert_{1} < \epsilon\right\}$,
fournissent également de bons $L$-modèles pour la partition plus grossière $\beta$
via
%
\begin{equation}
\left(\begin{array}{ccc}\mathcal{M}_\mu(\alpha, F, \epsilon, \sigma_n) & \longrightarrow & \mathcal{M}_\mu(\beta, F, \epsilon, \sigma_n) \\ a & \mapsto  & b=\Theta_{\beta,\alpha}\circ a\end{array}\right).
\end{equation}
Et ce sont ces modèles images qu'on dénombre. En d'autres termes, on considère le nombre de $L$-partitions de $D_n$ dans $\mathcal{M}_\mu(\beta, F, \epsilon, \sigma_n)$ qui sont susceptibles de se raffiner en des $K$-partitions de $D_n$ dans $\mathcal{M}_\mu(\alpha, F, \epsilon, \sigma_n)$, et on en étudie le taux de croissance exponentielle en la taille de $D_n$.

Cela conduit à une définition générale de l'entropie mesurée, sans hypothèse d'existence d'une partition génératrice finie.
Cette définition est due à D.~Kerr  \cite{Kerr-2013-Sofic-meas-ent-via-finite-partitions}. Elle est équivalente aux définitions de Kerr--Li \cite{Kerr-Li-2011-Variationnal-principle, Kerr-Li-2013-sofic-amenabl-dyn-entrop}, et toutes généralisent la définition due à L.~Bowen \cite{Bowen-2010-mes-conj-inv-sofic}.
\begin{defi}[Entropie sofique] 
L'\textbf{entropie sofique mesurée} de $G \curvearrowright (X, \mu)$ relativement à $\Sigma$ est définie comme
\begin{equation}
h_{\mathrm{mes}}^{\Sigma}(G\curvearrowright X,\mu) \overset{\mathrm{def}}{=} \sup_\beta \ \ \inf_{\alpha \geqslant \beta} \ \ \inf_{\epsilon > 0} \ \ \inf_{\substack{F \subset G\\F \text{ fini}}}\ \  \limsup_{n \rightarrow \infty} \frac{ \log \left| \Theta_{\beta,\alpha} \circ \mathcal{M}_\mu(\alpha, F, \epsilon, \sigma_n) \right|}{\vert D_n \vert},\tag{\textbf{Entropie sofique mesurée}}
\end{equation}
où  $\alpha$ et  $\beta$ parcourent les partitions mesurables finies de $X$.
\end{defi}
Force est d'admettre que cette formule est assez épouvantable. Mais en présence de partitions génératrices, on a de sérieuses simplifications. On donne quelques éléments supplémentaires en section~\ref{sect: ent sof. Kerr-Li th. generateur}.

Lorsque le groupe $G$ est moyennable infini, bien que ces définitions soient externes, utilisant des modèles finis, il se trouve qu'elles vont néanmoins coïncider avec les versions classiques de l'entropie \cite{Bowen-2012-sofic-ent-amenab, Kerr-Li-2013-sofic-amenabl-dyn-entrop}~:
\begin{theo}[Entropie pour les groupes moyennables]
\label{th: entropie sof. gp moyennable=KS}
Soit $G$ un groupe \textbf{moyennable} infini.
Pour toute action p.m.p. $G\curvearrowright^{\!T}{\!}(X, \mu)$, l'entropie mesurée classique de Kolmogorov-Sina{\u\i} et Ornstein-Weiss coïncide avec l'entropie mesurée sofique relativement à n'importe quelle approximation sofique $\Sigma$~:
\begin{equation}
h_{\mathrm{mes}}^{\Sigma}(G\curvearrowright X,\mu)=h_{\mathrm{KS}} (G\curvearrowright X,\mu).\tag{\textbf{Entropie des moyennables}}
\end{equation}
\end{theo}

\subsection{Entropie de Rokhlin}
\label{sect: ent Rokhlin}

L'origine de la notion d'\textbf{entropie de Rokhlin} est à chercher dans la majoration\footnote{Immédiate vu les propriétés de sous-additivité de la fonction $H$.}  $h_{\mathrm{KS}}(\mathbf{Z}\curvearrowright X,\mu)\leqslant H(\alpha)$, pour toute partition génératrice $\alpha$, et dans l'optimisation qu'en constitue le théorème des {\em générateurs de Rokhlin}\footnote{Concernant la manière d'orthographier son nom dans l'alphabet latin, observons que la plupart de ses articles en anglais sont publiés sous le nom de Rohlin. 
Il est cependant très plausible qu'il ait fini par préférer le nom de Rokhlin, sous lequel il a signé quelques-uns de ses derniers travaux, sous lequel il est référencé dans MathSciNet,  dans le  {\em Mathematics Genealogy Project} ou sur Wikipedia. C'est également cette orthographe qui est utilisée pour la traduction de sa notice nécrologique, dans les articles historiques rédigés par A.~Vershik, et par son propre fils Vladimir Rokhlin Jr, professeur en informatique à Yale.} 
(\cite{Rohlin-1963-generators} ou les notes~\cite{Rohlin-1967-lectures-entropy})
~:
\begin{theo}[des générateurs de Rokhlin]\label{th:gen Rokhlin}
Si $G\curvearrowright(X,\mu)$ est une action p.m.p. libre ergodique de $G=\mathbf{Z}$, alors son entropie de Kolmogorov-Sina{\u\i} est l'infimum des entropies de Shannon de ses partitions génératrices~:
\begin{equation}
h_{\mathrm{KS}}(G\curvearrowright X,\mu)=\inf\{H(\alpha): \textrm{ $\alpha$ partition génératrice dénombrable}\}.
\tag{\textbf{Rokhlin}}
\end{equation}
\end{theo}
En particulier, si l'entropie est finie, alors il doit exister une partition génératrice dénom\-brable d'entropie finie.
En fait, le théorème du {\em générateur fini} de Krieger \cite{Krieger-1970-finite-generator} affirme l'existence, dans ce contexte, d'une partition finie à $k$ pièces sitôt que $h_{\mathrm{KS}}(\mathbf{Z}\curvearrowright X,\mu)\leqslant \log k$.

Le théorème de Rokhlin a connu plusieurs généralisations à des actions d'autres groupes, 
au premier rang desquels les \textbf{groupes abéliens} par J.-P.~Conze \cite{Conze-1972}.
Ce n'est que récemment qu'il a été formellement étendu à tous les groupes $G$ infinis moyennables par B.~Seward et R.~Tucker-Drob  \cite{Seward-Tucker-Drob-2014-arxiv}.

Cela conduit à la définition suivante d'entropie qui pourrait potentiellement rendre les mêmes services que les notions vues précédemment, sans restriction aucune sur la nature du groupe dénombrable. Elle est introduite et étudiée dans une série d'articles extrêmement prometteurs de B.~Seward \cite{Seward-2014-Krieger-finite-th-Rokhlin-1, Seward-2015-Krieger-finite-th-Rokhlin-2}.
\begin{defi}[Entropie de Rokhlin {\cite{Seward-2014-Krieger-finite-th-Rokhlin-1}}]\label{def:Entropie de Rokhlin(Seward)}
L'\textbf{entropie de Rokhlin} d'une action ergodique p.m.p. $G\curvearrowright^{\! T}{\!} (X,\mu)$ d'un groupe dénombrable infini quelconque est
définie comme~:
\begin{equation}
h_{\mathrm{Rok}}(G\curvearrowright^{\! T}{\!} X,\mu)\overset{\mathrm{def}}{=}\inf\{H(\alpha): \textrm{ $\alpha$ partition génératrice dénombrable}\}.
\tag{\textbf{Entropie de Rokhlin}}
\end{equation}
\end{defi}
Une version non ergodique à l'étude \cite{Seward-2016-Krieger-finite-th-Rokhlin-3} fait appel à l'entropie de Shannon conditionnelle relativement à la sous-$\sigma$-algèbre des parties $G$-invariantes.

Observons que si $\alpha \colon X\to K$ est une partition finie à $\vert K \vert$ pièces, alors le cardinal de $\mathcal{M}_\mu(\alpha, F, \epsilon, \sigma_n)$ est trivialement majoré par le nombre $\vert K \vert^{\vert D_n\vert}$ de toutes les $K$-partitions sur $D_n$, ce qui conduit à la majoration  $h_{\mathrm{mes}}^{\Sigma}(G\curvearrowright X,\mu)\leqslant \log \vert K \vert$.
Une estimée plus précise découle de l'article de L.~Bowen \cite[Prop. 5.3]{Bowen-2010-mes-conj-inv-sofic}~:~pour toute action ergodique p.m.p. $G\curvearrowright^{\! T}{\!} (X,\mu)$ et toute approximation sofique $\Sigma$ de $G$,
on a
\begin{equation}
h_{\mathrm{mes}}^{\Sigma}(G\curvearrowright^{\! T}{\!} X,\mu)\leqslant h_{\mathrm{Rok}}(G\curvearrowright^{\! T}{\!} X,\mu).\tag{\textbf{Entropie sofique vs Rokhlin}}
\end{equation}
L'entropie de Rokhlin prend manifestement ses valeurs dans $[0, +\infty]$.
On ignore si, en dehors des cas où $h_{\mathrm{mes}}^{\Sigma}(G\curvearrowright X,\mu)=-\infty$, l'inégalité ci-dessus peut être remplacée par une égalité. 
C'est tout de même ce qui se produit pour les décalages de Bernoulli à base finie ou dénombrable des groupes sofiques, puisque pour la partition canonique $\alpha  \colon  x\mapsto x(\mathrm{id}_{G})$, on a 
$H(\alpha)= h_{\mathrm{mes}}^{\Sigma}(G\curvearrowright K^G,\nu^{\otimes G})\leqslant h_{\mathrm{Rok}}(G\curvearrowright K^G,\nu^{\otimes G})\leqslant H(\alpha)$.

Si une partition d'entropie finie réalise l'infimum dans la définition de l'entropie de Rokhlin, on a une réciproque.

\begin{theo}[Seward {\cite[Cor. 1.4]{Seward-2015-Krieger-finite-th-Rokhlin-2}}]
Soit $G\curvearrowright^{\! T}{\!} (X,\mu)$ une action p.m.p. libre ergodique d'un groupe infini dénombrable  et $\alpha$ une partition génératrice dénombrable. Si $h_{\mathrm{Rok}}(G\curvearrowright^{\! T}{\!} X,\mu)=H(\alpha)<\infty$, alors l'action est conjuguée à un\footnote{Elle est en fait conjuguée au décalage de Bernoulli évident.} décalage de Bernoulli.
\end{theo}
\begin{theo}[Seward {\cite{Seward-2015-Krieger-finite-th-Rokhlin-2}}]
Soit $G$ un groupe infini dénombrable.
Si $G$ admet des actions libres ergodiques d'entropie de Rokhlin finie, arbitrairement grande, alors~:
\begin{itemize}
\item [(i)] l'entropie de Rokhlin des décalages de Bernoulli de $G$ est égale à l'entropie de Shannon de leur base ;
\item[(ii)] les facteurs des Bernoulli  de $G$ sont d'entropie de Rokhlin non nulle ;
\item[(iii)] $G$ satisfait la conjecture de surjonctivité de Gottschalk\footnote{Voir section~\ref{sect:groupes sofiques}.}. 
\end{itemize}
\end{theo}
De plus, B. Seward montre \cite[Cor. 1.14]{Seward-2015-Krieger-finite-th-Rokhlin-2} que si tout groupe dénombrable admet une action libre ergodique avec $h_{\mathrm{Rok}}(G\curvearrowright X,\mu)>0$, alors les propriétés (i, ii, iii) sont vérifiées pour tout groupe infini dénombrable.

Après des améliorations quantitatives pour $\mathbf{Z}$, dues notamment à Denker \cite{Denker-1974} et Grillenberger et Krengel \cite{Grillenberger-Krengel-1976-Krieger-gen-th}, le théorème du générateur fini de Krieger a été poussé à degré de généralité optimale par B.~Seward (pour l'entropie de Rokhlin et un groupe non nécessairement sofique), qui montre une certaine flexibilité dans les partitions génératrices (et incidemment qu'on peut concocter une version de l'entropie de Rokhlin avec des partitions finies).
\begin{theo}[Seward {\cite{Seward-2014-Krieger-finite-th-Rokhlin-1}}]
Soit $G$ un groupe infini dénombrable et  $G\curvearrowright^{\! T}{\!} (X,\mu)$ une action p.m.p. ergodique mais pas nécessairement libre sur $(X,\mu)$ sans atome.
Pour tout vecteur de probabilité (fini ou infini) $\bar{p}=(p_i)_{i\in K}$ tel que  $\protect{h_{\mathrm{Rok}}(G\curvearrowright^{\! T}{\!} X,\mu)< H(\bar{p})}$, il existe une partition génératrice dont les pièces sont exactement de mesure $\mu(A_i)=p_i$ pour tout~$i\in K$.
\end{theo}

\subsection{Entropie topologique et principe variationnel}
\label{sect: ent topologique Kerr-Li}

L'entropie s'est invitée en dynamique topologique avec l'article \cite{Adler-Konheim-McAndrew-1965-top-entropy} de R.~Adler, A.~Konheim et M.~McAndrew qui introduisent l'entropie topologique pour les homéomorphismes $S \colon  X\to X$ d'un espace compact métrisable\footnote{Ils considèrent plus généralement une application continue d'un espace topologique, mais pour ce qui nous intéresse, ce degré de généralité n'est pas requis.}, guidés par la méthode de Kolmogorov-Sina{\u\i}. Il s'agit d'un invariant de conjugaison topologique.
Le rôle des partitions est joué par les recouvrements ouverts $\mathfrak{A}$ de $X$, auxquels on attache un nombre $N(\mathfrak{A})\overset{\mathrm{def}}{=}$ le cardinal minimal d'un sous-recouvrement. Ils définissent~:
\begin{align*}
h_{\mathrm{top}}(\mathbf{Z}\curvearrowright^{\! S}{\!} X; \mathfrak{A})& \overset{\mathrm{def}}{=} \lim\limits_{n\to \infty} \frac {1}{n}\log N(\mathfrak A\vee S^{-1}\mathfrak A\vee\cdots\vee S^{-n}\mathfrak A),
\\
h_{\mathrm{top}}(\mathbf{Z}\curvearrowright^{\! S}{\!} X)& \overset{\mathrm{def}}{=} \sup\limits_{\mathfrak{A}} \{h_{\mathrm{top}}(S\curvearrowright X; \mathfrak{A}) : \mathfrak{A}\textrm{ recouvrement ouvert} \}\tag{\textbf{Entropie topologique}}.
\end{align*}

Une variante est proposée par R.~Bowen\footnote{Rufus Bowen, apparemment sans lien de famille avec Lewis Bowen, l'inventeur de l'entropie sofique mesurée.} \cite{Bowen-1971-ent-top} et E.~Dinaburg \cite{Dinaburg-1970-announct-of-Dinaburg-1971, Dinaburg-1971-conn-var-ent-charac} ; et c'est plutôt celle-ci qui inspirera D.~Kerr et H.~Li. Il s'agit de compter le nombre de segments d'orbites qui sont \textbf{$\kappa$-séparés} via une distance auxiliaire $\rho$.
Plus précisément, soit $N_{\kappa}(n, \rho_{\infty})$ le nombre maximum (de points $x\in X$ et) de fonctions\footnote{\em Segments d'orbites.} \[\phi_{x} \colon \{0, 1, \cdots, n\}\to X, \ i\mapsto S^{i}(x)\] qui soient deux à deux à $\rho_{\infty}$-distance $\geqslant \kappa>0$
où
\[\rho_{\infty}(\phi_x, \phi_y) \overset{\mathrm{def}}{=} \max\limits_{i\in \{0,1, \cdots, n\}} \rho(S^{i}(x), S^{i}(y)).\]
R.~Bowen \cite{Bowen-1971-ent-top} et E.~Dinaburg \cite{Dinaburg-1971-conn-var-ent-charac} montrent que le taux de croissance exponentielle en $n$ permet de retrouver l'entropie topologique d'Adler-Konheim-McAndrew et ce,  indépendamment du choix d'une distance $\rho$ compatible avec la topologie~:
\begin{equation}
h_{\mathrm{top}}(\mathbf{Z}\curvearrowright^{\! S}{\!} X)=\lim\limits_{\kappa\to 0} \limsup\limits_{n\to \infty} \frac{\log N_{\kappa}(n, \rho_{\infty})}{n}.
\tag{\textbf{Bowen-Dinaburg}}
\end{equation}

D.~Kerr et H.~Li  étendent la notion aux actions continues des groupes sofiques, sur les compacts métrisables. Leur première version \cite{Kerr-Li-2011-Variationnal-principle}, exprimée en termes d'algèbres d'opérateurs, est exprimée dans \cite{Kerr-Li-2013-sofic-amenabl-dyn-entrop} en termes dynamiques sur l'espace.

Soit $G$ un groupe sofique et $\Sigma = (\sigma_n  \colon  G \rightarrow \mathrm{Sym}(D_n))_{n}$ une approximation sofique de $G$.
Soit $G \curvearrowright X$ une action continue de $G$ sur un espace compact métrisable $X$. 

Soit $\rho  \colon  X \times X \rightarrow [0, \infty[$ une  \textbf{pseudo-distance}\footnote{Cette généralisation est particulièrement pertinente lorsqu'on regarde un décalage de Bernoulli $K^G$ et une pseudo-distance induite par sa partition canonique $\alpha \colon K^G\to K$.}
 continue sur $X$ qui soit \textbf{génératrice} ; \emph{i.e.} $\rho$ est symétrique, satisfait l'inégalité triangulaire, et pour tout $x \neq y \in X$ il existe un $g \in G$ avec $\rho(g \cdot x, g \cdot y) > 0$. Pour des applications $\phi, \phi'  \colon  D_n \rightarrow X$ on définit
\begin{equation*}
\rho_2(\phi, \phi') \overset{\mathrm{def}}{=} \left( \frac{1}{\vert D_n \vert} \sum_{v\in D_n} \rho(\phi(v), \phi'(v))^2 \right)^{1 / 2}
\hskip10pt \textrm{ et } \hskip10pt
\rho_\infty(\phi, \phi') \overset{\mathrm{def}}{=} \max_{v \in D_n} \rho(\phi(v), \phi'(v)).
\end{equation*}
Pour une partie finie $F \subset G$ et $\delta > 0$, soit
\begin{equation*}
\mathrm{Map}(\rho, F, \delta, \sigma_n) \overset{\mathrm{def}}{=} \left\{\phi  \colon  D_n \rightarrow X:
\forall f\in F,\ \ \rho_2 \left( \phi \circ \sigma_n(f), \ f \cdot \phi \right) < \delta\right\}.
\end{equation*}
C'est la collection des applications  $\phi \in X^{D_n}$ qui sont  {\em presque équivariantes} 
(à $\delta$ près, sous la \og presque action\fg\ $\sigma_n$ restreinte à la partie finie $F\subset G$).
Finalement, on pose~:
\begin{equation}
N_\kappa(\mathrm{Map}(\rho, F, \delta, \sigma_n), \rho_\infty)
\end{equation}
le cardinal maximal d'un ensemble \textbf{$(\rho_\infty, \kappa)$-séparé} ; \emph{i.e.} 
un ensemble tel que pour toute paire d'éléments $\phi$ et $\phi'$ on ait $\rho_\infty(\phi, \phi') \geqslant \kappa$. 
\begin{defi}[Kerr-Li {\cite[Def. 2.3]{Kerr-Li-2013-sofic-amenabl-dyn-entrop}}]
L'\textbf{entropie sofique topologique} de l'action continue $G \curvearrowright X$ sur le compact métrisable $X$, relativement à l'approximation sofique $\Sigma$, est définie comme~:
\begin{equation}
h_{\mathrm{top}}^{\Sigma}(G\curvearrowright X) = \sup_{\kappa > 0} \inf_{\delta > 0} \inf_{\substack{F \subset G \\ F \text{ fini}}} \limsup_{n \rightarrow \infty} \ \frac{\log N_\kappa \left( \mathrm{Map}(\rho, F, \delta, \sigma_n), \rho_\infty \right)}{\vert D_n \vert}  .
\tag{\textbf{Entropie sofique topologique}}
\end{equation}
\end{defi}
D.~Kerr et H.~Li \cite{Kerr-Li-2013-sofic-amenabl-dyn-entrop} démontrent que la valeur de $h_{\mathrm{top}}^{\Sigma}(G\curvearrowright X)$ ne dépend pas du choix de la pseudo-distance continue génératrice $\rho$ (bien qu'elle puisse dépendre de $\Sigma$). Observons que de nouveau $h_{\mathrm{top}}^{\Sigma}(G\curvearrowright X)\geqslant 0$ ou bien $h_{\mathrm{top}}^{\Sigma}(G\curvearrowright X) = -\infty$.

Lorsque le groupe $G$ est moyennable, ils montrent qu'on retrouve la notion classique d'entropie topologique, et ce pour tout choix d'approximation sofique
\cite{Kerr-Li-2013-sofic-amenabl-dyn-entrop}. C'est la version topologique du théorème~\ref{th: entropie sof. gp moyennable=KS}.

\medskip
Ce qu'on appelle le \textbf{principe variationnel} est un énoncé qui affirme que l'entropie topologique d'une action continue $G \curvearrowright^{\! S}\! X$ est le \emph{supremum} des entropies mesurées pour tous les éléments de $M(G \curvearrowright^{\! S}\! X)$~: l'ensemble des mesures boréliennes de probabilité $G$-invariantes.
Sa version classique pour les actions continues de $\mathbf{Z}$ est due à 
T.~Goodman \cite{Goodman-1971-princ-variat} et repose sur des résultats de E.~Dinaburg \cite{Dinaburg-1970-announct-of-Dinaburg-1971, Dinaburg-1971-conn-var-ent-charac} et W. Goodwyn\footnote{Qui démontre l'inégalité $\geqslant$.} \cite{Goodwyn-1969-maj-ent-top-by-mes}.

\begin{theo}[Principe variationnel, Kerr-Li {\cite[Th. 6.1]{Kerr-Li-2011-Variationnal-principle}}]
Soit  $G \curvearrowright^{\! S}\! X$ une action continue, sur le compact métrisable $X$, du groupe sofique $G$ et soit $\Sigma$ une approximation  sofique de $G$ ; alors 
\begin{equation}
h_{\mathrm{top}}^{\Sigma}(G \curvearrowright^{\! S}\! X)=\sup\left\{h_{\mathrm{mes}}^{\Sigma}(G \curvearrowright^{\! S}\! X, \mu) : \mu \in M(G \curvearrowright^{\! S}\! X)\right\}.\tag{\textbf{Principe variationnel}}
\end{equation}
\end{theo}
Observons qu'il se peut que de telles actions n'aient pas de mesure invariante. Il s'agit alors d'un cas où $h_{\mathrm{top}}^{\Sigma}(G \curvearrowright^{\! S}\! X)=-\infty$. Cela n'arrive jamais si $G$ est moyennable et d'ailleurs l'existence d'une action continue telle que $M(G \curvearrowright^{\! S}\! X)=\emptyset$ est un critère de non-moyennabilité.

La recherche et l'identification de mesure qui réalise le \emph{supremum} dans le principe variationnel, et les liens avec le nombre de points périodiques constituent des thèmes récurrents en entropie topologique.
On fournit quelques éléments dans ce sens dans le cadre sofique.

Par exemple, pour un décalage de Bernoulli $K^G$ de base finie d'un groupe sofique, l'entropie sofique mesurée maximale $=\log\vert K\vert$ est réalisée par la mesure $\nu_{\mathbf{u}}^{\otimes G}$ provenant de la probabilité uniforme $\nu_{\mathbf{u}}$ sur $K$.
On peut observer qu'il s'agit de la mesure de Haar sur $\mathbf{K}^G$ lorsque $K=\mathbf{K}$ est un groupe fini.

\begin{theo}[Gaboriau-Seward {\cite[Th. 8.2]{Gab-Seward-2015-arxiv}}]
\label{th:Gab-Sew haar-ent max}
Si $H$ est un groupe profini sur lequel le groupe sofique $G$ agit par automorphismes continus, de sorte que le sous-groupe homocline soit dense, alors la mesure de Haar de $H$ est d'entropie maximale~:
\[h_{\mathrm{top}}^{\Sigma}(G\curvearrowright H)=h_{\mathrm{mes}}^{\Sigma}(G\curvearrowright X, \mathrm{Haar}),\]
pour toute approximation sofique $\Sigma$ de $G$.
\end{theo}
On rappelle que le sous-groupe \textbf{homocline} est le sous-groupe des points $h\in H$ tels que $g_n.h\to \mathrm{id}_{H}$ pour toute suite injective $(g_n)_{n}$ dans $G$, c'est par exemple les éléments de support fini dans $H<\mathbf{K}^{G}$.

\begin{theo}[Gaboriau-Seward {\cite[Th. 4.6]{Gab-Seward-2015-arxiv}}]
\label{th: Gab-Sew ent et pt fixes}
Soit $G$ un groupe résiduellement fini et $\Sigma$ une approximation sofique associée à une chaîne $(G_n)_{n}$ de sous-groupes d'indice fini. Soit $\mathbf{K}$ un groupe fini et $X\subset \mathbf{K}^G$ un sous-groupe compact $G$-invariant 
qui soit un sous-décalage de type fini.
Alors 
\[h_{\mathrm{top}}^{\Sigma}(G\curvearrowright X)=\limsup \frac{1}{[G:G_n]} \ \log \vert \mathrm{Fix}_{G_n} (X) \vert. \]
\end{theo}
Rappelons qu'une partie  $X\subset \mathbf{K}^G$ est un \textbf{sous-décalage de type fini} (défini par une partie finie $W$ du groupe $G$ et une partie $P\subset \mathbf{K}^W$) 
si c'est un fermé $G$-invariant qui est maximal sous la condition
que pour tout $x\in X$, l'application $W\to \mathbf{K}$, $w\mapsto \alpha(w.x)$ appartienne à $P$.

Si $\mathbf{K}$ et $\mathbf{L}$ sont des groupes finis et $\Phi \colon \mathbf{K}^G\to \mathbf{L}^G$ est un homomorphisme $G$-équi\-variant continu, alors on peut appliquer 
le Théorème~\ref{th:Gab-Sew haar-ent max} à l'image de $\Phi$ et le Théorème~\ref{th: Gab-Sew ent et pt fixes} au noyau de $\Phi$.

\section{Définitions \& notations}
\label{sect: definitions}

Si $Y$ est un ensemble fini, on note $\vert Y\vert $ son cardinal.
Pour un ensemble fini non vide $D$, on note $\mathbf{u}$ la \textbf{mesure uniforme} $\mathbf{u}(A)=\frac{\vert A\vert}{\vert D\vert}$ sur $D$.
On désigne par $\mathrm{Sym}(Y)$ le groupe symétrique (de toutes les permutations) de $Y$.

\subsection{Partitions}

Une \textbf{partition} $\alpha=(A_k)_{k\in K}$ d'un ensemble $X$ est une famille 
 de parties de $X$ mutuellement disjointes et qui forme un recouvrement de $X$.
En particulier, une permutation des indices conduit en général à des partitions distinctes.
Pour insister sur ce point, on utilise parfois la terminologie \textbf{$K$-partition}.

Un autre point de vue intéressant consiste à considérer la partition $\alpha$ de $X$ comme la fonction $\alpha \colon  X\to K$ qui à $x\in X$ associe l'indice de la pièce qui le contient, et donc $A_k=\alpha^{-1}(k)$.
Dans ce contexte, $\alpha$ est parfois appelée une \textbf{observable}, notamment chez L.~Bowen \cite{Bowen-2010-mes-conj-inv-sofic}.

Dans le cadre d'un espace de probabilité $(X,\mathcal{B}, \mu)$, c'est-à-dire un ensemble $X$ muni d'une tribu $\mathcal{B}$ (qu'on omet de mentionner lorsqu'une confusion nous paraît improbable) 
et d'une mesure de probabilité $\mu$ sur $\mathcal{B}$, on se restreint à des partitions au plus dénombrables, dont les pièces sont mesurables.

L'\textbf{entropie de Shannon} d'une partition finie ou dénombrable $\alpha=(A_k)_{k\in K}$ est définie par la formule~:
\begin{equation}
H(\alpha)\overset{\mathrm{def}}{=}-\sum_{k\in K} \mu(A_k)\log \mu(A_k).
\end{equation}
Elle représente, en théorie de l'information, la quantité d'information ($-\log \mu(A_k)$ pour chaque pièce $A_k$) contenue en moyenne dans les pièces de la partition.
Le sort des pièces de mesure nulle est réglé en convenant que $0 \log(0)=0$.

Si $\nu$ est une mesure de probabilité sur un ensemble fini ou dénombrable $K$, on note encore $H(\nu)$
l'entropie de la partition en singletons.

On note $\alpha\vee \beta \colon 
\left(\begin{array}{ccc} 
X & \to & K\times L
\\
x&\mapsto &(\alpha(x), \beta(x))
\end{array}\right)$ le \textbf{joint} des partitions $\alpha \colon X\to K$ et $\beta \colon  X\to L$, c'est-à-dire la partition formée des intersections
$A_{k}\cap B_{l}=\alpha^{-1}(k)\cap \beta^{-1}(l)$.

La fonction $t\mapsto -t\log t$ pour $t\in ]0,1]$ est concave. Lorsqu'on raffine une partition, son entropie augmente.
On a toujours, 
$H((A_k)_{k\in K})\leqslant \log \vert K\vert $ avec égalité si et seulement si les parties $A_k$ 
ont même mesure.
Par ailleurs, $H(\alpha\vee \beta)\leqslant H(\alpha)+H(\beta)$ avec égalité si et seulement si les partitions sont indépendantes~:  $\mu(A_i\cap B_j)=\mu(A_i)\mu(B_j)$ pour tout $(i,j)\in K\times L$.

\subsection{Partitions et actions de groupes}

De manière générale, une \textbf{conjugaison} entre deux actions $G\curvearrowright^{\! T} \! X$ et $G\curvearrowright^{\! S} \! Y$ est un isomorphisme $\Phi \colon  X\to Y$ 
tel que 
\begin{equation}
\tag{\'Equivariance}
 g\in G,  \ \ \Phi\circ T(g)(x)=S(g)\circ \Phi(x).
\end{equation}
Les actions sont alors dites \textbf{conjuguées}.

Si on parle d'actions continues, on demande que $\Phi$ soit un homéomorphisme.

Si on parle d'action p.m.p. $G\curvearrowright^{\! T} \! (X, \mu)$ et $G\curvearrowright^{\! S} \! (Y, \nu)$, on demande une bijection bimesurable préservant la mesure entre deux parties $X'\subset X$ et $Y'\subset Y$ de mesure pleine telle que la condition d'équivariance soit vérifiée pour tout $x\in X'$.

Dans le cadre mesuré, un \textbf{facteur} est une application mesurable équivariante $\Phi \colon  (X,\mu)\to (Y,\nu)$ essentiellement surjective~: $Y\setminus \Phi(X)$ est négligeable.
%
\subsection{Décalages de Bernoulli}

De façon générale, si $K$ est un ensemble et $V$ un ensemble dénombrable muni d'une action $G\curvearrowright V$ d'un groupe dénombrable, alors l'espace $K^{V}=\prod_{v\in V} K$ des fonctions $V\to K$ est muni de l'action par décalage~:
\[\forall x\in K^G, \forall g\in G, \ g\cdot x(v)= x(g^{-1} v), \ \forall v\in V.\]
Lorsque $V=G$ sur lequel $G$ agit par multiplication à gauche, alors $G\curvearrowright K^{G}$ est appelé \textbf{décalage de Bernoulli} de base $K$.

Si $K$ est un espace topologique métrisable séparable, $K^{V}$ est équipé de la topologie produit et l'action par décalage est continue.

Une mesure borélienne $\nu$ sur $K$ délivre la mesure borélienne produit $\nu^{\otimes G}$ sur $K^G$. Elle est invariante sous l'action de $G$. On parle encore de \textbf{décalage de Bernoulli} (mesuré, cette fois-ci)~:
\[  G\curvearrowright (K^{G}, \nu^{\otimes G}).\]
Lorsque $K$ est fini ou dénombrable muni d'une mesure $\nu$, la partition \textbf{canonique} $\alpha$ est définie par l'évaluation 
$\alpha \colon  \left(\begin{array}{ccc} 
K^{G}&\to &K\\
  x&\mapsto &x(\mathrm{id}_{G})
\end{array}\right)$ en l'élément neutre du groupe.
 Elle est génératrice au sens ci-dessous.

\subsection{Partitions génératrices}
\label{subsubsect:partitions generatrices}
Une action p.m.p. $G\curvearrowright (X,\mathcal{B}_X,\mu)$ étant donnée, une $K$-partition finie ou dénombrable $\alpha \colon X\to K$ nous parle en réalité d'un \textbf{facteur} de cette action. 
Elle nous fournit la sous-$\sigma$-algèbre $\mathcal{S}_{\alpha}$ (engendrée par les $\cup_{F\subset G, F \textrm{ finie }} \alpha^{\vee F}$)
laissée globalement invariante par $G$ et à laquelle correspond un certain facteur.
On peut  introduire explicitement le facteur~: c'est l'application $G$-équivariante naturelle $\pi \colon X\to K^{G},  x\mapsto (\alpha(g^{-1}(x)))_{g\in G}$, où $G$ agit sur $K^G$ par décalage de Bernoulli et laisse invariante la mesure $\nu\overset{\mathrm{def}}{=}\pi_{*}(\mu)$, image directe de la mesure $\mu$.

Une partition $\alpha$ est \textbf{génératrice} \index{génératrice! partition}\index{partition! génératrice} s'il existe une partie $X'\subset X$ de mesure $1$ telle que pour  tout $ x\not= y\in X'$, il y a un $g\in G$ pour lequel $\alpha$ sépare $g.x$ de $g.y$ \emph{i.e.} $\alpha(g.x)\not=\alpha(g.y)$.
De manière équivalente, la partition $\alpha$ est génératrice si et seulement si $\pi$ est essentiellement injectif (\emph{i.e.} injectif sur une partie de mesure pleine).

Inversement, pour un facteur $\pi \colon (X,\mathcal{B}_X,\mu)\to (Y,\mathcal{B}_Y,\nu)$ non essentiellement injectif, toute partition finie mesurée $\xi$ de $Y$ se relève en une partition $\alpha=\pi^{-1}(\xi)$ non génératrice~: la sous-$\sigma$-algèbre engendrée n'est pas capable de séparer les points d'une même fibre.

\section{Groupes sofiques}\label{sect:groupes sofiques}

Un groupe dénombrable est résiduellement fini s'il admet une \textbf{chaîne normale}, c'est-à-dire une suite décroissante de sous-groupes normaux d'indices finis $(G_n)_{n}$ telle que $\cap_n G_n=\{1_G\}$.
Ce qu'on va en retenir, c'est qu'on dispose d'ensembles finis ${D_n\overset{\mathrm{def}}{=}G_n\backslash G}$ (disons les classes à droite) et d'applications dans leurs groupes symétriques $\sigma_n \colon  G \rightarrow \mathrm{Sym}(D_n)$ (obtenues par multiplication à droite par l'inverse\footnote{Voir ci-dessous l'interprétation en termes de graphes, pour une explication sur ce choix.}), telles que
\begin{enumerate}
\item[\rm (i)] (action) les $\sigma_n$ sont des homomorphismes de groupes,
\item[\rm (ii)] (liberté) les $\sigma_n$ séparent les éléments de $G$  en un sens fort~: pour tout $g\in G\setminus\{1_G\}$ et $n$ assez grand, l'image $\sigma_n(g)$ agit sans point fixe sur l'ensemble fini $D_n$.
\end{enumerate}
Ce sont ces conditions qu'on imite, en les relaxant, dans la définition~\ref{defn: gp approx sofic} pour obtenir la notion de groupe sofique, en leur demandant d'être satisfaites asymptotiquement.
On demande que la proportion de points de $D_n$ sur lesquels elles sont vérifiées tende vers~$1$.

Si on préfère une version non asymptotique, on dira que $G$ est \textbf{sofique} si 
pour toute partie finie $F\subset G$ et tout réel $\delta>0$, il existe une \textbf{$(F,\delta)$-approximation}, c'est-à-dire un ensemble fini $D$ et une application $\sigma \colon G \rightarrow \mathrm{Sym}(D)$ tels que 
\begin{enumerate}
\item[\rm (i)] (presque action) \hfil
$\frac{1}{\vert D\vert}\, \vert \{ {v} \in D \ \vert \ \sigma_n(g) \circ \sigma_n(h)({v}) = \sigma_n(g h)({v})\} \vert \geqslant  1-\delta$, \hfil $\forall g, h \in F$
\item[\rm (ii)] (presque libre) \hfil
$\frac{1}{\vert D\vert}\, \vert\{ {v} \in D \ \vert \ \sigma_n(g)({v}) \neq {v}\}\vert \geqslant  1-\delta$, \hfil $\forall g\in F\setminus\{1_G\}$.
\end{enumerate}

On constate que la condition de normalité de la chaîne de sous-groupes d'indices finis n'est pas indispensable pour que les $G_n\backslash G$ fournissent une approximation sofique. La condition optimale (parfois appelée \textbf{condition de Farber}) est que $\forall g \in G \setminus \{1_G\}$, 
$ \lim_{n \rightarrow \infty} \mathbf{u}_{G_n\backslash G} \left(\{G_n v \ \vert \ v \in G, \textrm{ t. q. } \ g \in v^{-1} G_n v\}\right)= 0$.
Dans ce cas, on dira que $(G_n)_{n}$ est une  \textbf{chaîne sofique}.

Les groupes résiduellement finis, par exemple les groupes linéaires de type fini (Malcev 1940) et notamment les groupes libres, sont  sofiques.

On dispose d'une autre grande classe de groupes sofiques~: les groupes moyennables.
Rappelons qu'un groupe dénombrable est \textbf{moyennable} s'il admet une \textbf{suite de F{\o}lner}, c'est-à-dire une suite de parties finies $F_n\subset G$ qui vérifient 
pour tout $g\in G$
\begin{equation*}
\lim\limits_{n\to \infty}\frac{\vert F_n\cdot g^{-1}\Delta F_n\vert}{\vert F_n\vert}=0.
\end{equation*}
Autrement dit, l'action de $G$ sur lui-même, par multiplication à droite par l'inverse, laisse les parties $F_n$ asymptotiquement invariantes.
On obtient alors une approximation sofique de $G$ en observant que, par multiplication à droite par l'inverse, tout élément $g$ de $G$ définit une bijection de l'ensemble fini $D_n\overset{\mathrm{def}}{=}F_n\subset G$, à une partie asymptotiquement négligeable près. On pose alors $\sigma_n(g)(f)=f g^{-1}$ lorsque $f$ et $f g^{-1}$ appartiennent à $F_n$ et on l'étend à souhait en une bijection de $F_n$ tout entier.

Voici une liste de quelques propriétés de stabilité pour la soficité.
Un groupe est sofique si et seulement si tous ses sous-groupes de type fini sont sofiques.
Si un groupe $G$ possède un sous-groupe normal sofique tel que le quotient soit moyennable, alors $G$ est sofique.
Les produits directs de groupes sofiques sont sofiques.
Un produit amalgamé ou une HNN-extension de groupes sofiques au-dessus d'un groupe moyennable est sofique (G.~Elek et E.~Szab{\'o} \cite{Elek-Szabo-2011-sofic-amalg-over-amenable} et  indépendamment L.~ P{\u{a}}unescu \cite{Paunescu-2011-sofic-act-equiv-rel}).

Une interprétation agréable, notamment lorsque le groupe $G$ est engendré par une partie finie $S$, consiste à considérer $D_n$ comme l'ensemble des sommets d'un graphe $\mathcal{G}_n$ dont les arêtes, orientées et étiquetées par $s\in S$, joignent chaque $v$ à son image par $\sigma_n(s)$~:
\begin{equation}
\mathcal{G}_{n}=(D_n, ([v, \sigma_n(s).v])_{v\in D_n, s\in S}).
\end{equation}
Les conditions (i) et (ii) de la définition de soficité reviennent à dire que les $v\in D_n$ depuis lesquels la boule de centre $v$, de rayon $R$, dans $\mathcal{G}_n$ est isomorphe\footnote{Comme graphe orienté, étiqueté.} à la boule de même rayon\footnote{Disons centrée en $\mathrm{id}_{G}$, mais les boules de même rayon sont toutes isomorphes.} dans le graphe de Cayley\footnote{Où $\rho$ est l'action à gauche de $G$ sur lui-même par multiplication à droite par l'inverse. Un \textbf{graphe de Cayley} est un graphe orienté, équipé d'une action de $G$ simplement transitive sur les sommets. Le choix d'un point base permet d'identifier les sommets avec $G$ et fournit une action $\rho$ qui commute avec la première, et dont découle l'étiquetage.}
\begin{equation}
\mathrm{Cayley}(G, S, \rho)=\mathcal{G}=(G, ([v, \rho(s).v])_{g\in G, s\in S})
\tag{\textbf{Graphe de Cayley}} 
\end{equation}
forment une partie de $D_n$ dont, à $R$ fixé, la proportion  tend vers $1$ lorsque $n$ tend vers l'infini~:
\begin{equation}
\forall R>0, \ \ \lim\limits_{n \to \infty}\mathbf{u}_{n}\left(\left\{ v\in D_n\ \vert\  B_{\mathcal{G}_{n}}(v,R)\simeq B_{\mathcal{G}}(v,R)\right\}\right)=1.
\tag{\textbf{Soficité \& graphes}}
\end{equation}
Les exemples décrits ci-dessus s'interprètent alors de la manière suivante.
\\
-- Si $G$ est résiduellement fini et $(G_n)_{n}$ est une chaîne sofique de sous-groupes, alors $\mathcal{G}_{n}$ n'est autre que ce qu'on appelle le graphe de Schreier, le graphe quotient~:
\begin{equation}
\mathcal{G}_{n}=G_n\backslash \mathrm{Cayley}(G, S, \rho).
\tag{\textbf{Graphe de Schreier}} 
\end{equation}
-- Si $G$ est moyennable, alors $\mathcal{G}_{n}$ est formé à partir de la restriction du graphe  $\mathrm{Cayley}(G, S, \rho)$ à la partie de F\o lner $F_n\subset G$, en bricolant les arêtes du bord, c'est-à-dire en reliant bijectivement pour chaque $s\in S$ les sommets $\{g\in F_n \ \vert\ gs^{-1}\not\in F_n\}$ aux sommets $\{gs^{-1} \in F_n \ \vert\ g\not\in F_n\}$.
\\
-- Les graphes associés aux groupes libres, pour $S$ partie génératrice libre, sont ceux dont le tour de taille tend vers l'infini (voir section~\ref{sect:plus sur l'ex de OW}).

\medskip
Les groupes sofiques ont été introduits (sous un autre nom\footnote{Groupes à graphes de Cayley {\em initially subamenable}.}) par M.~Gromov \cite{Gromov-1999-symbolic-algebraic-varieties},
pour lesquels il a montré la validité de la \textbf{conjecture de surjonctivité de Gottschalk}~: toute application continue $G$-équivariante injective $\Phi \colon K^G\to K^G$ est automatiquement surjective, où $K$ est fini et $G$ agit par décalage de Bernoulli.

D.~Kerr et H.~Li en donnent une preuve entropique.
%
\begin{theo}[Kerr-Li {\cite[Th. 4.2]{Kerr-Li-2011-Variationnal-principle}}] 
\label{th: Kerr-Li entrop Gottschalk conject}
Soit $G \curvearrowright K^{G}$ l'action (continue) d'un groupe sofique $G$ par décalage de Bernoulli de base finie.
Toute restriction de cette action à une partie fermée propre $G$-invariante, est d'entropie sofique topologique strictement inférieure à celle de $G\curvearrowright K^G$.
\end{theo}
En particulier, si $\Phi \colon X=K^{G}\to K^{G}$ est une application continue $G$-équivariante injective, son image est un fermé $G$-invariant de même entropie que $G\curvearrowright K^{G}$. 
Cette image ne peut pas être propre ; $\Phi$~doit être surjective.

\medskip
La terminologie \textbf{sofique}, dérivée d'un mot hébreu signifiant \og fini\fg, a été introduite par B.~Weiss \cite{Weiss-2000-sofic-gp}.
Il faut noter qu'à ce jour, on ne connaît aucun exemple de groupe qui ne soit pas sofique.
Pour une jolie introduction aux groupes sofiques, on pourra consulter \cite{Pestov-2008-sofic-gps-survey}. Voir aussi l'ouvrage \cite{Capraro-Lupini-LNM-sofic-hyperlin-gps}.

\section{Entropie sofique, le point de vue externe}

\subsection{Un exemple de \og point de vue externe\fg\ sur l'entropie, sans dynamique}
\label{sect: preuve du lem Boltzmann-Sanov}

On revient sur la proposition~\ref{prop:Boltzmann-Sanov} \og Boltzmann-Sanov\fg\ et l'entropie sans action de groupe.
 
Soit $\alpha=(A_k)_{k\in K}$ une $K$-partition finie d'un espace de probabilité $(X,\mu)$. 
On considère, sur un ensemble fini\footnote{Qu'on pense grand ; et ce qu'on regarde ne dépend que de son cardinal.} $D$ de cardinal $d$, la famille de toutes les partitions $(V_k)_{k\in K}$  qui imitent bien $\alpha$ en termes de mesure des pièces pour la mesure uniforme sur $D$
\begin{equation}
 {\mathcal{M}(\alpha,\epsilon,d)}\overset{\mathrm{def}}{=}\left\{K\textrm{-partitions } (V_k)_{k\in K} \textrm{ de } D \ \Big\vert \ \forall k\in K,\  \left\vert \frac{\vert V_k \vert}{\vert D\vert} -\mu(A_k)\right\vert <\epsilon\right\},   \label{eq: estimee exp nb bon. part.}
\end{equation}
pour un certain $\epsilon>0$ petit. On compare le taux de croissance exponentielle en $d$ de son cardinal à l'entropie $H(\alpha)$. 
Plus précisément~:
\begin{lemm}[Entropie de Shannon et modèles finis, Boltzmann-Sanov]
\label{lem:entropie sans groupe}
Dans ce contexte, 
$\forall\delta>0$,  
$\exists\epsilon_0>0$ tel que $\forall \epsilon \in ]0,\epsilon_0[$,
$\exists d_0\in \mathbf{N}$ tel que pour tout ensemble fini $D$ de cardinal $d>d_0$~:
\begin{equation}
\mathrm{e}^{ (H(\alpha)-\delta) d}
<  
\left\vert
\mathcal{M}(\alpha,\epsilon,d)
\right\vert
< \mathrm{e}^{  (H(\alpha)+\delta) d}. 
\end{equation}
\end{lemm}
Et cela conduit à une définition \og externe\fg\ de l'entropie de Shannon de $\alpha$, en extrayant le taux de croissance exponentielle du nombre de bons modèles (Proposition~\ref{prop:Boltzmann-Sanov}).

On a des preuves purement combinatoires de ce lemme~\ref{lem:entropie sans groupe}, en utilisant la formule de Stirling et les estimées standard.
Mais, en supposant que les $\mu(A_k)$ sont tous non nuls, on a aussi une preuve, qui dans un esprit de principe de grande déviation adopte une mesure adéquate, et qui s'avérera utile en section~\ref{sect: ent sofique Bernoulli}.

Au fond,  en interprétant les $K$-partitions de $D$ comme des fonctions $a \colon D\to K$,
\begin{equation*}
\mathcal{M}(\alpha,\epsilon,d)=\left\{
a\in K^{D} \ \Big\vert \ \forall k\in K, \   \left\vert \frac{\vert a^{-1}(k)\vert}{\vert D\vert}  - \mu(A_k) \right\vert < \epsilon 
\right\}
\end{equation*}
et
cette mesure $\frac{\vert a^{-1}(k)\vert}{\vert D\vert} =\frac{\vert \{v\in D \ \vert \ a(v)=k\} \vert}{\vert D\vert}$ n'est autre que la fréquence d'apparition de la lettre $k$ dans la \og suite\fg\ $(a(v))_{v\in D}$. Et ce qu'on requiert dans (\ref{eq: estimee exp nb bon. part.}), c'est que cette fréquence soit proche de $\mu(A_k)$. Or on dispose d'une mesure naturelle sur $K^D$ qui est bien adaptée à cette question, via la loi des grands nombres~:
quand on effectue un grand nombre de tirages aléatoires indépendants dans $(K,\nu_{\alpha})$, où la mesure  $\nu_{\alpha}$ sur $K$ est définie par $\nu_{\alpha}(\{k\})=\mu(A_k)\overset{\mathrm{def}}{=}p_k$, une grande proportion de ces tirages verra la fréquence d'apparition de la lettre $k$ proche de $\mu(A_k)$. En d'autres termes, avec la mesure produit ${\nu_{\alpha}}^{\otimes D} $ sur $K^D$, on a~: 
$\forall \delta>0, \forall\epsilon>0$,
 il existe une taille $d_1$ de $D$ à partir de laquelle~:
\begin{align*}
{\nu_{\alpha}}^{\otimes D} \left(\left\{a\in K^D\ \Bigg\vert \ 
{ \forall k
\textrm{ la fréquence de $k$ dans $(a(v))_{v\in D}$ }   
\atop
\textrm{    est dans $]p_k-\epsilon, p_k+\epsilon[$}
}\right\}
\right) 
>
 1 -\delta.
\end{align*}
Mais maintenant qu'on en connaît à peu près les fréquences, chaque atome de notre ensemble $\mathcal{M}(\alpha,\epsilon,d)$ a une mesure qui vérifie~:
\begin{equation}
\prod_{k\in K} p_k^{(p_k+\epsilon) d} <  {\nu_{\alpha}}^{\otimes D} (\textrm{atome de } \mathcal{M}(\alpha,\epsilon,d)) <  \prod_{k\in K } p_k^{(p_k-\epsilon) d} ;
\end{equation}
soit, $\forall \delta>0$, dès que $0<\epsilon< \frac{\delta}{-2 \sum_{k\in K} p_k}$,  pour tout $d\geqslant d_1$ l'encadrement suivant~:
\begin{equation}\label{eq: mesure d'un atome de Bernoulli Dn}
\mathrm{e}^{(-H(\alpha) -\delta/2)d}
 <  {\nu_{\alpha}}^{\otimes D} (\textrm{atome de } \mathcal{M}(\alpha,\epsilon,d)) <  
\mathrm{e}^{(-H(\alpha)+\delta/2)d}.
\end{equation}
Ce qui nous délivre les estimations  sur le cardinal~:
\begin{equation}
\mathrm{e}^{(-H(\alpha) -\delta/2)d} \vert \mathcal{M}(\alpha,\epsilon,d)\vert 
<  {\nu_{\alpha}}^{\otimes D} ( \mathcal{M}(\alpha,\epsilon,d)) 
< \mathrm{e}^{(-H(\alpha) +\delta/2)d}
 \vert \mathcal{M}(\alpha,\epsilon,d)\vert.
\end{equation}
Et pour $d$ suffisamment grand~:
\begin{equation}
\mathrm{e}^{(H(\alpha) -\delta)d}< (1-\delta) 
\mathrm{e}^{(H(\alpha) -\delta/2)d}
<  \vert \mathcal{M}(\alpha,\epsilon,d)\vert 
< 
\mathrm{e}^{(H(\alpha) +\delta/2)d}
< 
\mathrm{e}^{(H(\alpha) +\delta)d}.
\end{equation}

\subsection{Les modèles finis de la dynamique}
\label{sect: modeles finis de la dynamique}
On ajoute l'ingrédient d'une approximation sofique $\Sigma$ d'un groupe sofique $G$ et
une action p.m.p. $G\curvearrowright (X,\mu)$.

Au risque d'être redondant,  une partition finie $\alpha  \colon  X \rightarrow K$ et la partie finie $F \subset G$ délivrent pour chaque  $x\in X$ un certain élément $p\in K^{F}$ ;  le \textbf{$F$-parcours} $f\mapsto \alpha(f\cdot x)$ de $x$.

De façon analogue, une partition $a  \colon  D_n \rightarrow K$ délivre pour tout  ${v}\in D_n$ un certain élément $p\in K^{F}$ ;
le $F$-(pseudo)-parcours $f\mapsto a(\sigma_n(f)({v}))$ de $v$.

Ces données définissent ainsi des partitions de $X$ (resp. $D_n$), selon le parcours associé, en les pièces suivantes indexées par les $p\in K^{F}$~:
\begin{align}
U_{p}(\alpha,F)&\overset{\mathrm{def}}{=}\{x\in X\ \vert\ \forall f\in F, \ \alpha(f\cdot x)=p(f)\},\label{def:U p alpha F}
\\
U_{p}(a,F,n)&\overset{\mathrm{def}}{=}\{{v}\in D_n\ \vert\ \forall f\in F, \ a(\sigma_n(f) \cdot {v})=p(f)\} \label{def:U p a F n}.
\end{align}
Et on compare les mesures $\mu(U_{p}(\alpha,F))$ et $\mathbf{u}_n(U_{p}(a,F,n))$ de ces pièces. Plus précisément, pour tout  $\epsilon > 0$, posons~:
\begin{equation}\label{eq:def M mu}
\mathcal{M}_\mu(\alpha, F, \epsilon, \sigma_n)\overset{\mathrm{def}}{=}\left\{a \in K^{D_n}\ \Big\vert \ \sum_{p \in K^F}
\left| \mu(U_{p}(\alpha,F)) - \mathbf{u}_n(U_{p}(a,F,n)) 
\right| < \epsilon\right\},
\end{equation}
l'ensemble de toutes les partitions $a\in K^D$ qui sont  $(F,\epsilon)$-plausibles dans l'approximation $\sigma_{n} \colon G\to \mathrm{Sym}(D_{n})$ de $G$. Et on s'intéresse à leur nombre $\vert \mathcal{M}_\mu(\alpha, F, \epsilon, \sigma_n)\vert $.

On a des propriétés de monotonie immédiates~:
\begin{equation}
\textrm{si $F'\subset F$ et $\epsilon\leqslant \epsilon'$, alors }
 \mathcal{M}_\mu(\alpha, F, \epsilon, \sigma_n)\subset \mathcal{M}_\mu(\alpha, F', \epsilon', \sigma_n).
\end{equation}

\subsection{Entropie sofique sans partition génératrice}
\label{sect: ent sof. Kerr-Li th. generateur}

On rappelle que si $\alpha  \colon  X \rightarrow K$ est une partition plus fine que 
$\beta  \colon  X \rightarrow L$, l'application de fusion des pièces est notée $\Theta_{\beta,\alpha}  \colon  K \rightarrow L$
(voir section~\ref{subsect: Ent sof sans part gen}).

\begin{defi}[Entropie sofique mesurée]\label{def:ent sof mes Kerr}
Soit $\Sigma$ une approximation sofique du groupe $G$ et soit $G\curvearrowright(X, \mathcal{B}_X,\mu)$ une action action p.m.p.
L'\textbf{entropie sofique mesurée} de l'action relativement à $\Sigma$ est définie comme~:
\begin{equation*}
h_{\mathrm{mes}}^{\Sigma}(G\curvearrowright X, \mu)\overset{\mathrm{def}}{=}\sup_{\beta}\  
\underbrace{\inf_{\alpha\geqslant \beta}\  
\underbrace{\inf_{F} \inf_{\epsilon>0} 
\underbrace{\limsup_{n\to \infty} \frac{\log \vert \Theta_{\beta,\alpha} \circ \mathcal{M}_\mu(\alpha, F, \epsilon, \sigma_n) \vert}{\vert D_n\vert}}_{h_{\mathrm{mes}}^{\Sigma}(\beta, \alpha, F, \epsilon)}}_{h_{\mathrm{mes}}^{\Sigma}(\beta, \alpha)}}_{h_{\mathrm{mes}}^{\Sigma}(\beta)},
\end{equation*}
où  $\alpha$ et  $\beta$ parcourent les partitions mesurables finies et $\alpha$ est plus fine que $\beta$.
\end{defi}
\begin{theo}[Générateur]
\label{th: generateur -Kerr}
Soit $\Sigma$ une approximation sofique du groupe $G$ et soit $G\curvearrowright(X, \mathcal{B}_X,\mu)$ une action action p.m.p.
Si $\beta$ est une partition génératrice finie, alors elle réalise le \emph{supremum}~:
\begin{equation}
h_{\mathrm{mes}}^{\Sigma}(G\curvearrowright X, \mu)=h_{\mathrm{mes}}^{\Sigma}(\beta).
\end{equation}
De plus, pour toute partition finie plus fine $\alpha\geqslant \beta$, on a~:
\begin{equation}
h_{\mathrm{mes}}^{\Sigma}(\beta)=h_{\mathrm{mes}}^{\Sigma}(\beta,\alpha)=h_{\mathrm{mes}}^{\Sigma}(\beta, \beta).
\end{equation}
\end{theo}
\begin{coro}
Si $\alpha$ est une partition génératrice, alors on retrouve bien la définition de L.~Bowen~:
\begin{equation*}
h_{\mathrm{mes}}^{\Sigma}(G\curvearrowright X,\mu) =  \inf_{\epsilon > 0} \inf_{\substack{F \subset G\\F \text{ fini}}} \limsup_{n \rightarrow \infty} \frac{\log \left|  \mathcal{M}_\mu(\alpha, F, \epsilon, \sigma_n) \right|}{\vert D_n \vert} .
\end{equation*}
\end{coro}
Le théorème~\ref{th: generateur -Kerr} fait usage des quatre ingrédients suivants, dont les preuves reposent sur des arguments de dénombrement sans énorme surprise, mais assez délicats.
Les entités mises en jeu sont des partitions finies $\alpha, \beta, \xi$ telles que $\alpha\geqslant \beta$ et $W$ est une partie finie du groupe, qui contient $\mathrm{id}_{G}$. 
\begin{enumerate}
\item 
\label{ingredient iterer beta}
$h_{\mathrm{mes}}^{\Sigma}(\beta)=h_{\mathrm{mes}}^{\Sigma}(\beta^{\vee W})$. 
\item 
\label{ingredient : ineg de Rokhlin}
$h_{\mathrm{mes}}^{\Sigma}(\xi)\leqslant h_{\mathrm{mes}}^{\Sigma}(\beta) +H(\xi \vert \beta)$  pour toute $\xi$ (\og Inégalité de Rokhlin\fg). 
\item $h_{\mathrm{mes}}^{\Sigma}(\beta,\alpha^{\vee W})=h_{\mathrm{mes}}^{\Sigma}(\beta, \alpha)$. 
\label{ingredient iterer alpha}
\item 
\label{ingredient partitions dans sous-algebre}
Soit $\mathcal{S}$  une sous-algèbre dense dans $\mathcal{B}_X$ qui contient $\beta$. Pour toute $\alpha\geqslant \beta$, il existe une partition $\xi$ dans $\mathcal{S}$ telle que $\xi\geqslant \beta$
telle que $h_{\mathrm{mes}}^{\Sigma}(\beta, \xi)\leqslant h_{\mathrm{mes}}^{\Sigma}(\beta, \alpha)$.
\end{enumerate}
Les ingrédients \ref{ingredient iterer beta} et \ref{ingredient iterer alpha} étudient le comportement des quantités lorsque les partitions sont itérées sous $W$.
Les deux autres permettront des arguments d'approximations de partitions quelconques par des partitions de la forme $\beta^{\vee W}$.
Le point~\ref{ingredient partitions dans sous-algebre} notamment est rassurant puisqu'il autorisera à effectuer les calculs dans une sous-algèbre de parties.
\\
Les ingrédients \ref{ingredient iterer beta} et \ref{ingredient : ineg de Rokhlin} suffisent\footnote{Et c'est désormais le même genre d'argument que pour le théorème de Kolmogorov-Sina{\u\i}.}
 à montrer que les partitions finies génératrices réalisent l'entropie sofique (le \emph{supremum} des $h_{\mathrm{mes}}^{\Sigma}(\xi)$). 
En effet si $\beta$ est génératrice, alors pour toute partition $\xi$, l'entropie relative $H(\xi \vert \beta^{\vee W})$ peut être rendue aussi petite que l'on souhaite en choisissant $W$ partie finie assez grande de $G$. Alors, 
\begin{equation*}
h_{\mathrm{mes}}^{\Sigma}(\xi)\leqslant h_{\mathrm{mes}}^{\Sigma}(\beta^{\vee W}) +H(\xi \vert \beta^{\vee W})=h_{\mathrm{mes}}^{\Sigma}(\beta) +H(\xi \vert \beta^{\vee W})\underset{''W\nearrow G''}{\longrightarrow}h_{\mathrm{mes}}^{\Sigma}(\beta).
\end{equation*}
Les ingrédients \ref{ingredient iterer alpha} et \ref{ingredient partitions dans sous-algebre} fournissent les outils qui permettront les calculs.
Si $\beta$ est génératrice, l'algèbre $\cup_{\substack{W\subset G \\ W \textrm{ finie}}} \beta^{\vee W}$ peut jouer le rôle de $\mathcal{S}$ dans \ref{ingredient partitions dans sous-algebre}.
Ainsi, pour toute $\alpha\geqslant \beta$, on peut trouver $W$ et $\beta^{\vee W} \geqslant \xi\geqslant \beta$ tels que $h_{\mathrm{mes}}^{\Sigma}(\beta, \xi) \leqslant h_{\mathrm{mes}}^{\Sigma}(\beta, \alpha)$ ; et alors
\begin{equation*}
\begin{array}{ccccccc} 
h_{\mathrm{mes}}^{\Sigma}(\beta,\beta)  \overset{(\ref{ingredient iterer alpha})}{=}h_{\mathrm{mes}}^{\Sigma}(\beta,\beta^{\vee W}) & \overset{(*)}{\leqslant} &  h_{\mathrm{mes}}^{\Sigma}(\beta, \xi) &\overset{(\ref{ingredient partitions dans sous-algebre})}{\leqslant}& h_{\mathrm{mes}}^{\Sigma}(\beta, \alpha) & \overset{(*)}{\leqslant}  & h_{\mathrm{mes}}^{\Sigma}(\beta,\beta), \end{array}
\end{equation*}
où les deux inégalités $(*)$ reposent sur la monotonie évidente~: si $\alpha'\geqslant \alpha$, alors $h_{\mathrm{mes}}^{\Sigma}(\beta,\alpha')\leqslant h_{\mathrm{mes}}^{\Sigma}(\beta, \alpha)$.

\subsection{Actions Bernoulli}\label{sect: ent sofique Bernoulli}
On donne des éléments de preuve permettant le calcul de l'entropie sofique mesurée des décalages de Bernoulli.
\begin{theo}[Bowen {\cite{Bowen-2010-mes-conj-inv-sofic}}]
Soit $G$ un groupe sofique infini et $\Sigma$  une approximation sofique de $G$.
Soit $K$ un ensemble fini, muni d'une mesure de probabilité $\nu$ qui charge chacun de ses points 
et $G\curvearrowright (K^{G}, \nu^{\otimes G})$ l'action par décalage de Bernoulli associée. Alors,
\begin{equation}
h_{\mathrm{mes}}^{\Sigma}\left(G\curvearrowright K^{G}, \nu^{\otimes G}\right)=H(\nu).
\end{equation}
\end{theo}
Considérons la partition canonique $\alpha \colon K^{G}\to K$,  $x\mapsto x(\mathrm{id}_{G})$.
Elle est génératrice et son entropie de Shannon vaut
$H(\alpha)=H(\nu)$.

On se restreint au cas où le groupe $G$ est résiduellement fini et où l'approximation sofique vérifie~:
\\
-a- les $\sigma_n \colon G\to \mathrm{Sym}(D_n)$ sont des homomorphismes, et 
\\
-b- pour toute partie finie $F\subset G$,  si $n$ est assez grand, alors \textbf{pour tout} $v\in D_n$ et tout $s,t\in F$, si $s\not=t$ alors 
$\sigma_n(s).v\not= \sigma_n(t).v$.

D'après le Théorème \&\ Définition~\ref{th-def: ent sof mes Bowen}, on cherche à estimer le cardinal de l'ensemble (voir (\ref{def:U p alpha F}), (\ref{def:U p a F n}) et (\ref{eq:def M mu}))~:
\begin{equation}
\mathcal{M}_\mu(\alpha, F, \epsilon, \sigma_n)\overset{\mathrm{def}}{=}\left\{a \in K^{D_n}\ \Bigg\vert \ \sum_{p \in K^F}
\left| \mu\left(U_{p}(\alpha,F)\right) - \mathbf{u}_n\left(U_{p}(a,F,n)\right) 
\right| < \epsilon\right\}.
\end{equation}

On le regarde comme partie de l'espace probabilisé  $(K^{D_n}, \nu^{\otimes D_{n}})$.
On sait déjà que, pour cette mesure déjà considérée $\nu^{\otimes D_{n}}$, quand $\mathrm{id}_{G}\in F$ chaque atome est de mesure approximativement $\mathrm{e}^{-H(\nu) \vert D_n\vert}$ (voir section~\ref{sect: preuve du lem Boltzmann-Sanov}).
On va montrer, par une variante standard de la loi des grands nombres (avec dépendance limitée),
 que
pour toute partie finie $F\subset G$ et pour tout $\epsilon>0$, avec une très grande probabilité lorsque $D_n$ est grand, tous les $a \in K^{D_n}$ vont convenir.
Plus précisément :
pour tout $p\in K^{F}$, pour tout $\epsilon>0$, 
\begin{equation}
\nu^{\otimes D_{n}}\left(\left\{a \in K^{D_n}\ \vert\  
\left| \mu(U_{p}(\alpha,F)) - \mathbf{u}_n(U_{p}(a,F,n)) 
\right|< \epsilon\right\}\right)\underset{n\to \infty}{\longrightarrow} 1.
\end{equation}
Fixons $p\in K^{F}$. 
On cherche à estimer, pour chaque $a\in K^{D_n}$, le cardinal
\begin{equation}
\left\vert U_{p}(a,F,n)\right\vert =\left\vert \left\{{v}\in D_n\ \vert\ \ a^{\vee F}(v)=p\right\}\right\vert =\sum_{v\in D_n}
\mathbf{1}_{a^{\vee F}(v)=p},
\end{equation}
 qu'on regarde comme variable aléatoire sur $(K^{D_n}, \nu^{\otimes D_{n}})$
\begin{equation}
Z \colon a\mapsto \sum_{v\in D_n} Z_{v}(a)
\ \ \textrm{ en posant } \ \  Z_{v}(a)\overset{\mathrm{def}}{=} \mathbf{1}_{a^{\vee F}(v)=p}.
\end{equation}
On va utiliser l'inégalité de Bienaymé-Tchebychev\footnote{Si $V$ une variable aléatoire d'espérance $m$ et de variance finie $\sigma^2$, 
alors pour tout réel strictement positif $\delta$, on a~: $P\left(\left|V-m\right| \geqslant \delta \right) \leqslant \frac{\sigma^2}{\delta^2}$.}
qui majore la probabilité de l'écart à la moyenne à l'aide de la variance.

\noindent\textbf{Espérance de $Z$.}
\begin{equation}
\mathbb{E}_{a}(Z)= \mathbb{E}_{a}\left(\sum_{v\in D_n} Z_{v}(a)\right)=\sum_{v\in D_n} \mathbb{E}_{a}\left(Z_{v}(a)\right)=\vert D_n\vert \ \mu\left(U_p(\alpha, F)\right).
\end{equation}
En effet, pour $n$ suffisamment grand, pour chaque $v\in D_n$, les images $\sigma_{n}(f).v$ pour $f\in F$ sont deux à deux distinctes. Du coup $\mathbb{E}_{a}(Z_{v}(a))$, qui est la mesure pour $\nu^{\otimes D_n}$ du cylindre de $K^{D_n}$ dont les valeurs aux coordonnées $\sigma_n(f_i).v$ sont fixées à $p(f_i)$, vaut exactement $\mu(U_p(\alpha, F))$, \emph{i.e.} la mesure d'un cylindre analogue pour $\nu^{\otimes G}$.

\medskip
\noindent\textbf{Variance de $Z$.}
Si $v,w\in D_n$ sont deux points dont les $F$-parcours sont disjoints $\{\sigma_n(f) \cdot v: f\in F\}\cap \{\sigma_n(f) \cdot w: f\in F\}=\emptyset$, alors les cylindres qu'ils peuvent définir portant sur des coordonnées distinctes, 
les variables aléatoires $Z_v$ et $Z_{w}$ sont indépendantes~: $\mathbb{E}_{a}(Z_{v}Z_{w})=\mathbb{E}_{a}(Z_{v})\ \mathbb{E}_{a}(Z_{w})$.
Le nombre de paires $(v,w)\in D_n^2$ pour lesquelles $Z_v$ et $Z_{w}$ ne sont pas indépendantes est donc inférieur à $\vert D_n\vert \vert F\vert^2$.
Ainsi puisque chacun des $ \mathbb{E}_{a}(Z_{v}\ Z_{w})\in [0,1]$, on a~:
\[\mathbb{E}_{a}(Z^2)=\sum_{(v,w)\in D_n} \mathbb{E}_{a}(Z_{v}\ Z_{w})\leqslant 
\sum_{(v,w)\in D_n} \mathbb{E}_{a}(Z_{v})\ \mathbb{E}_{a}(Z_{w})+ \vert D_n\vert \vert F\vert^2.\]
Cela donne une majoration de la variance de~$Z$~:
\[\mathrm{Var}(Z)= \mathbb{E}_{a}(Z^2)- \mathbb{E}_{a}(Z)^2\leqslant \vert D_n\vert \vert F\vert^2.\]
L'inégalité de Bienaymé-Tchebychev donne alors~:
\begin{align*}
\mathbb{P}_{a} \left( \left\vert  \frac{Z}{\vert D_n\vert} -\frac{\mathbb{E}_{a}(Z)}{\vert D_n\vert} \right\vert \geqslant \epsilon_1 \right) &\leqslant \frac{\vert F\vert ^2}{\epsilon_1 ^2 \vert D_n\vert}, \textrm{ \emph{i.e.} }
\\
\nu^{\otimes D_{n}}\left(\left\{a \in K^{D_n}\ \vert\  
\left| \mu\left(U_{p}(\alpha,F)\right) - \mathbf{u}_n\left(U_{p}(a,F,n)\right) 
\right|< \epsilon_1
\right\}\right) 
&> 1-\frac{\vert F\vert ^2}{\epsilon_1 ^2 \vert D_n\vert}
\end{align*}
D'où avec $\epsilon_1=\frac{\epsilon}{\vert K\vert^{\vert F\vert}}$
\[
\nu^{\otimes D_{n}}\left(\mathcal{M}_\mu(\alpha, F, \epsilon, \sigma_n)\right) > 1- \frac{\vert F\vert ^2}{\epsilon ^2 \vert D_n\vert}.
\]
Alors exactement comme en section~\ref{sect: preuve du lem Boltzmann-Sanov}, l'estimée sur la mesure de $\mathcal{M}_\mu(\alpha, F, \epsilon, \sigma_n)$ et celle des atomes permet de conclure.
Pour tout $\delta>0$, il existe $\epsilon_0$ tel que pour tout $\epsilon\in]0, \epsilon_0[$, pour toute partie finie $F\subset G$ avec $\mathrm{id}_{G}\in F$, il existe $n_0$ tel que pour tout $n\geqslant n_0$ on a~:
\begin{equation}
\mathrm{e}^{(H(\alpha) -\delta)d}
<  \vert  \mathcal{M}_\mu(\alpha, F, \epsilon, \sigma_n)\vert 
< \mathrm{e}^{(H(\alpha) +\delta)d}.
\end{equation}
Et donc
\begin{equation}
H(\nu)=H(\alpha)=\inf\limits_{\epsilon>0}\ \inf\limits_{\substack{F\subset G\\ \textrm{ fini}}}\ \limsup\limits_{n\to \infty} \frac{ \log\vert \mathcal{M}_\mu(\alpha, F, \epsilon, \sigma_n)\vert }{\vert D_n\vert}. 
\end{equation}

\section{Dépendance en l'approximation sofique}\label{sect: depend-approx-sofic}

On montre des situations où la valeur de l'entropie sofique dépend de l'approximation sofique.

Soit $G$ un groupe $G$ avec la \textbf{propriété (T) de Kazhdan}. C'est-à-dire qu'il vérifie~: 
$\forall \delta>0$, il existe une partie finie $F\subset G$ et $\epsilon>0$ tels que pour toute représentation unitaire $\pi$, si $\xi$ est un vecteur unitaire $(F,\epsilon)$-invariant (\emph{i.e.} $\Vert \pi(f) \xi-\xi\Vert ^2<\epsilon$), alors il existe un vecteur unitaire $G$-invariant $\xi_0$ tel que $\Vert \xi_0-\xi\Vert^2<\delta$.

Si $G$ a la propriété (T) de Kazhdan et est résiduellement fini, alors  une approximation sofique $\Sigma$ associée à une chaîne $(G_n)_{n}$ de sous-groupes normaux d'indices finis n'est pas du tout encline à modéliser des actions non ergodiques de $G$.

Plus précisément, si $G\curvearrowright(X,\mu)$ se décompose en deux parties $G$-invariantes $X=X_1\sqcup X_2$ de mesures non nulles, alors pour toute partition $\alpha$ qui raffine cette décomposition, pour $\epsilon$ assez petit et $F$ assez grand on a $ \mathcal{M}_\mu(\alpha, F, \epsilon, \sigma_n)=\emptyset$, pour tout $n$, et donc~:
\[h_{\mathrm{mes}}^{\Sigma}(G\curvearrowright X,\mu)=-\infty.\]
Cela repose sur le fait suivant~: considérons une action transitive $\sigma \colon  G\curvearrowright D$ sur un ensemble fini (par exemple sur $G_n\backslash G$). Pour $\delta>0$, si $A\subset D$ est  $(F,\epsilon)$-invariante (pour $(F,\epsilon)$ donnés par la propriété de Kazhdan) au sens où pour tout $f^{-1}\in F$
\begin{equation*}
\vert f^{-1}A\Delta A\vert <\epsilon \vert A\vert,
\end{equation*}
alors $A$ occupe une grande proportion de $D$~:
\begin{equation*}
\vert A\vert > (1-\delta/2)^2 \ \vert D\vert.
\end{equation*}
En effet, la représentation unitaire associée sur $\ell^2(D)$ définie pour $g\in G, \xi \in \ell^2(D)$ et $v\in D$ par $(\pi(g).\xi )(v)=\xi(\sigma(g). v)$ possède deux vecteurs unitaires $G$-invariants~:~$\pm \frac{\mathbf{1}_{D}}{\sqrt{\vert D\vert}}$.
Puisque $\Vert \mathbf{1}_{A}-\mathbf{1}_{B}\Vert^2=\vert A\Delta B\vert$ pour toutes parties $A,B\subset D$, on~a~:
$
\vert f^{-1}A\Delta A\vert =\Vert \mathbf{1}_{A}-\mathbf{1}_{f^{-1}(A)}\Vert^2=\Vert \mathbf{1}_{A}-\pi(f^{-1}).\mathbf{1}_{A}\Vert^2<\epsilon \Vert \mathbf{1}_{A}\Vert^2 =\epsilon \vert A\vert.
$ 
La propriété (T) montre alors que $\Vert  \frac{\mathbf{1}_{D}}{\sqrt{\vert D\vert}}-\frac{\mathbf{1}_{A}}{\sqrt{\vert A\vert}} \Vert^2<\delta$ soit $\frac{\sqrt{\vert A\vert}}{\sqrt{\vert D\vert}}=
\langle \frac{\mathbf{1}_{D}}{\sqrt{\vert D\vert}},\frac{\mathbf{1}_{A}}{\sqrt{\vert A\vert}} \rangle
> 1-\delta/2$.

 Admettons maintenant que les restrictions de notre action $G\curvearrowright (X,\mu)$ à $X_1$ et $X_2$ sont ergodiques (par exemple des  décalages de Bernoulli de base finie) et disons que $(\mu(X_1), \mu(X_2))=(2/5,3/5)$.
L'approximation sofique $\Sigma'$ consistant en cinq copies de $\Sigma$, \emph{i.e.} $D'_{n}=D_{n}\times\{1,2,3, 4, 5\}$ et $\sigma_{n}' \colon G\to \mathrm{Sym}(D_{n}\times\{1,2,3, 4, 5\})$ est induite par $\sigma_{n}$ sur la première coordonnée, sera elle bien adaptée.
Pour $\epsilon$ assez petit et $F$ assez grand, les bons modèles pour $G\curvearrowright (X_1\sqcup X_2,\mu_{1}+\mu_{2})$ vont se décomposer en bons modèles pour $G\curvearrowright (X_1,\mu_1)$ (sur deux des $D_{n}\times \{i\}$) et pour 
 $G\curvearrowright (X_2,\mu_2)$ (sur les trois autres $D_{n}\times \{j\}$). Leur nombre sera approximativement~:
 \[ C_5^2\,\left(\mathrm{e}^{\vert D_n\vert \,h_{\mathrm{mes}}^{\Sigma} (G\curvearrowright X_1, \mu_1)}\right)^{2} \left(\mathrm{e}^{\vert D_n\vert \,h_{\mathrm{mes}}^{\Sigma} (G\curvearrowright X_2, \mu_2)}\right)^{3}.\]
  Et finalement pour cette deuxième approximation sofique $\Sigma'$, on a~:
 \begin{eqnarray*}
h_{\mathrm{mes}}^{\Sigma'}\left(G\curvearrowright X, \mu\right)=\frac{2}{5} h_{\mathrm{mes}}^{\Sigma} \left(G\curvearrowright X_1, \mu_1\right) + \frac{3}{5} h_{\mathrm{mes}}^{\Sigma} \left(G\curvearrowright X_2, \mu_2\right),
\end{eqnarray*}
dont on peut facilement prescrire des valeurs finies $h_{\mathrm{mes}}^{\Sigma} \left(G\curvearrowright X_i, \mu_i\right)=t_i$ pour des décalages de Bernoulli.

Observons qu'on n'a pas utilisé toute la force de la propriété (T), mais seulement la propriété $(\tau)$
(voir le livre de A.~Lubotzky \cite{Lubotzky=book-exp-graph-1994}) du groupe $G$ vis-à-vis de sa suite de sous-groupes d'indice fini $(G_n)_{n}$ ; plus précisément que la représentation unitaire
 $G\curvearrowright \oplus_n \ell^{2}_{0}(G_n\backslash G)$ ne contient pas faiblement la représentation triviale, où $\ell^{2}_{0}(G_n\backslash G)$ est l'orthogonal des fonctions constantes.
En d'autres termes, on a utilisé le caractère \textbf{expanseur} des graphes de Schreier $\mathcal{G}_n$ de la section~\ref{sect:groupes sofiques}.

Le théorème suivant de A.~Carderi est frappant dans ce sens qu'il met en évidence une forme de rigidité dans l'adéquation chaîne sofique/action.
\begin{theo}[Carderi {\cite[Th. D]{Carderi-2015-arxiv}}]
Soit $G$ un groupe libre ou bien $\mathrm{PSL}_{r}(\mathbf{Z})$ pour $r\geqslant 2$.
Il existe un continuum de chaînes normales $(H_n^{t})_{n\in \mathbf{N}}$ de sous-groupes de $G$ (indexées par $t\in \mathbf{R}$) telles que les entropies sofiques mesurées et topologiques, relativement à l'approximation sofique $\Sigma_t$ associée à la chaîne $(H_n^{t})_{n\in \mathbf{N}}$, de l'action profinie $G\curvearrowright \left(\varprojlim(H_n^{s})_{n}, \mu_s\right)$ associée à la chaîne $(H_n^{s})_{n\in \mathbf{N}}$ (avec son unique mesure invariante $\mu_s$) vérifient 
\[
h_{\mathrm{top}}^{\Sigma_t}(G\curvearrowright \varprojlim(H_n^{s})_{n})=h_{\mathrm{mes}}^{\Sigma_t}(G\curvearrowright \varprojlim(H_n^{s})_{n},\mu_s)=\begin{cases} 0 & \textrm{ si } t=s,\\ -\infty & \textrm{ si } t\not=s. \end{cases}\]
\end{theo}
Le résultat repose via \cite{Abert-Elek=profinite-act-2012} sur la propriété $(T)$ de Kazhdan ou sur le fait que les sous-groupes de congruence $N_i$ de $\mathrm{PSL}_{2}(\mathbf{Z})$ ont la propriété $(\tau)$.
Les familles de sous-groupes $H_n=N_{i_1}\cap N_{i_2}\cap\cdots\cap N_{i_n}$ associées à des suites infinies incomparables (au sens de l'inclusion) d'indices $I=\{i_1, i_2, i_3, \cdots\}\subset \mathbf{N}$ conviendront. 

Notons que pour les actions profinies générales $G \curvearrowright \varprojlim(G_n)_{n}$, l'entropie sofique topologique (ou mesurée pour l'unique mesure invariante) prend ses valeurs dans $\{-\infty, 0\}$.

\section{Miscellanées}

\subsection{Un peu plus sur l'exemple d'Ornstein-Weiss}
\label{sect:plus sur l'ex de OW}
Soit  $\mathbf{L}_r$ le groupe libre à $r$ générateurs engendré par $S=\{a_1, a_2, \cdots, a_r\}$ et soit $\mathbb{K}$ un corps fini.
Soit $\mathcal{G}$ le graphe (l'arbre) de Cayley associé. L'application
\begin{equation*}
\Theta_{r}  \colon  \left(\begin{array}{ccl}
\mathbb{K}^{\mathbf{L}_r}&\to & (\mathbb{K}^{r})^{\mathbf{L}_r}\\
\omega &\mapsto &
 (\omega(g a_1)-\omega(g), \omega(g a_2)-\omega(g), \cdots, \omega(g a_r)-\omega(g))_{g\in \mathbf{L}_r}
\end{array}\right)
\end{equation*}
de la section~\ref{sect: appl Ornstein-Weiss}  (\ref{eq: appl OW}) s'interprète comme l'application linéaire cobord \[\delta^{1} \colon C^{0}(\mathcal{G}, \mathbb{K})\to C^{1}(\mathcal{G}, \mathbb{K})\] entre les cochaînes en dimension $0$ et les cochaînes en dimension $1$ à coefficients dans $\mathbb{K}$.
On constate facilement qu'elle est surjective.
L'entropie sofique topologique (ou mesurée pour les mesures de Haar) vaut $\log \vert \mathbb{K}\vert$ pour $\mathbb{K}^{\mathbf{L}_r}$, resp. $r\,\log \vert \mathbb{K}\vert$ pour $\mathrm{im}\, \delta^{1}$.

Les approximations sofiques $\Sigma=(\sigma_n \colon \mathbf{L}_r\to \mathrm{Sym}(D_n))_{n}$ du groupe libre $\mathbf{L}_r$ correspondent à des graphes $\mathcal{G}_{n}=(D_n, ([v, \sigma_n(s).v])_{v\in D_n, s\in S})$ dont le tour de taille (longueur minimale de cycles) tend vers l'infini avec $n$. Chaque sommet est de valence $2r$ et est le sommet initial et le sommet terminal d'une arête étiquetée $a_i$ pour chaque $i=1, 2, \cdots, r$.
On peut considérer l'application cobord 
\begin{equation}
\label{eq: appl cobord}
\begin{array}{rccc}
\delta^{1}_{n} \colon &C^{0}(\mathcal{G}_n, \mathbb{K})&\to& C^{1}(\mathcal{G}_n, \mathbb{K})
\\
&\simeq \mathbb{K}^{\vert D_n \vert} & & \simeq \mathbb{K}^{r\,\vert D_n \vert}
\end{array}
\end{equation}
 sur ces graphes.
Des considérations de dimensions, avec $\dim_{\mathbb{K}}C^{0}(\mathcal{G}_n, \mathbb{K})=\vert D_n\vert$, $\dim_{\mathbb{K}}C^{1}(\mathcal{G}_n, \mathbb{K})=r\,\vert D_n\vert$ et $\dim_{\mathbb{K}}\ker \delta^{1}_{n}=1$, 
nous indiquent par le théorème du rang ($\dim_{\mathbb{K}}\mathrm{im}\, \delta^{1}_{n}= \vert D_n\vert -1$) que l'image est loin d'être surjective !

On a vu que dans $C^{0}(\mathcal{G}_n, \mathbb{K})\simeq \mathbb{K}^{\vert D_n \vert}$, la plupart des points (au sens de la mesure de Haar) fournissent de bons modèles pour le calcul de l'entropie sofique de $G\curvearrowright \mathbb{K}^G$.
En revanche, leurs images par $\delta^{1}_{n}$ sont en nombre insuffisant pour représenter tous les bons modèles de $G\curvearrowright (\mathbb{K}^r)^{G}$.

Si l'approximation sofique est donnée par une chaîne de sous-groupes normaux $(G_n)_{n}$ d'indice fini, on peut tester le Théorème~\ref{th: Gab-Sew ent et pt fixes}~:
les points fixes de $G_n\curvearrowright C^{0}(\mathcal{G}, \mathbb{K})$ et de $G_n\curvearrowright C^{1}(\mathcal{G}, \mathbb{K})$ sont précisément les relevés de 
$C^{0}(G_n\backslash\mathcal{G}, \mathbb{K})$ et de $C^{1}(G_n\backslash\mathcal{G}_n, \mathbb{K})$ tandis que $\delta^{1}\left(\mathrm{Fix}_{G_n}C^{0}(\mathcal{G}, \mathbb{K})\right)$ est le relevé de
$\mathrm{im}\,\delta^{1}_{n}$~: bien trop \emph{petit}.

Observons que ces applications $\delta^{1}, \delta^{1}_{n}$ de (\ref{eq: appl cobord}) se généralisent à tout groupe avec $r$ générateurs\footnote{Mais $\delta^{1}$ n'est pas surjective en général.} pour le graphe de Cayley $\mathcal{G}$ associé, et qu'en vérité, l'action image s'identifie à 
\[G\curvearrowright C^{0}(\mathcal{G}, \mathbb{K})/\ker \delta^{1}\simeq \mathbb{K}^{G}/\mathbb{K}\] et ne dépend donc pas du système générateur. 
Le \og défaut de commutation dans le diagramme~: pousser par $\delta$, puis prendre les modèles finis ou prendre d'abord les modèles finis, puis pousser par $\delta$\fg\
a été exploité dans \cite[Th. 9.4]{Gab-Seward-2015-arxiv} pour interpréter la croissance de l'entropie sofique par le facteur $\delta^{1} \colon  \mathbb{K}^{G}\to \mathbb{K}^{G}/\mathbb{K}$ en termes de \emph{coût\footnote{On renvoie à \cite{Gab00a} pour cette notion.} du groupe $G$} ou en termes de $\beta_{(2)}^{1}(G)$, son \emph{premier nombre de Betti} $\ell^2$, via le théorème d'approximation de L\"uck et sa généralisation aux approximations sofiques \cite{Luc94b,Thom=Diophantine-approx=2008}. En particulier, lorsque $G$ est de type fini, $\delta_1$ fait croître l'entropie sitôt que le premier nombre de Betti $\ell^2$ de $G$ est non nul~:
\[(1+\beta_{(2)}^{1}(G))\, \log\vert \mathbb{K}\vert \leqslant h_{\mathrm{mes}}^{\Sigma}(G\curvearrowright \mathbb{K}^{G}/\mathbb{K})=h_{\mathrm{top}}^{\Sigma}(G\curvearrowright \mathbb{K}^{G}/\mathbb{K}).\]

La question du caractère Bernoulli\footnote{Est-elle conjuguée à un décalage de Bernoulli ?} de l'image est très largement ouverte.
Lorsque $G$ est moyennable, les facteurs des Bernoulli sont des Bernoulli \cite{OW87}. En revanche, des travaux de S.~Popa et R.~Sasyk \cite{PS07,Popa-1-cohomology-2006} montrent que pour un groupe infini avec la propriété (T), l'action $G\curvearrowright \mathbb{K}^{G}/\mathbb{K}$ n'est pas Bernoulli. Cela passe par le calcul explicite du premier groupe de cohomologie de ces actions (égal au groupe fini $\mathrm{Char}(G)$ des caractères de $G$  pour Bernoulli~; isomorphe à $\mathrm{Char}(G)\times \mathbb{K}$ pour le quotient).

\subsection{Produits}

L'entropie de Kolmogorov-Sina{\u\i} est additive sous produits cartésiens. 
Considérons deux actions p.m.p. $G \curvearrowright^{\! T}\! (X,\mu)$ et $G \curvearrowright^{\! S}\! (Y,\nu)$.
Si $G=\mathbf{Z}$, alors 
$$h_{\mathrm{KS}}(G \curvearrowright^{\! T\times S}\! X\times Y,\mu\times \nu)=h_{\mathrm{KS}}(G \curvearrowright^{\! T}\! X,\mu) + h_{\mathrm{KS}}(G \curvearrowright^{\! S}\! Y,\nu).$$
{\c C}a n'est plus vrai pour l'entropie sofique. Tim Austin \cite{Austin-2015-Add-prod-sofic} a montré l'inégalité
$$h_{\mathrm{mes}}^{\Sigma}(G \curvearrowright^{\! T\times S}\! X\times Y,\mu\times \nu)\leqslant h_{\mathrm{mes}}^{\Sigma}(G \curvearrowright^{\!T}\! X,\mu) + h_{\mathrm{mes}}^{\Sigma}(G \curvearrowright^{\! S}\! Y,\nu),$$
 et donné des contre-exemples à l'égalité.
Cependant, si l'un des deux facteurs est un décalage de Bernoulli, d'entropie de Shannon de base finie, alors on a égalité 
\cite[Th. 8.1]{Bowen-2010-mes-conj-inv-sofic}.

\subsection{Actions non libres}

Soit $G\curvearrowright (X,\mu)$ une action p.m.p. ergodique d'un groupe sofique et $\Sigma$ une approximation sofique.
Si $h_{\mathrm{mes}}^{\Sigma}(G\curvearrowright X, \mu)>0$, alors le stabilisateur de $\mu$-presque tout point est fini  \cite{Meyerovitch-2015-ent-stabilizer}.
L'entropie de Rokhlin n'a aucune prise sur ce genre de question, puisqu'elle ne voit pas les stabilisateurs.

\subsection{Actions algébriques et déterminant de Fuglede-Kadison}
 Un autre thème récurrent en théorie classique de l'entropie consiste à relier l'entropie de systèmes dynamiques d'origine algébrique aux valeurs propres d'une l'application linéaire sous-jacente.
Par exemple, pour un automorphisme $\phi$ linéaire hyperbolique\footnote{Les valeurs propres sont de module $\not=1$.} du tore $\mathbf{R}^{n}/\mathbf{Z}^{n}$ 
 $$h_{\mathrm{top}}(G\curvearrowright X_{\phi})=h_{\mathrm{KS}}(G\curvearrowright X_{\phi}, \mathrm{Haar})=\log{ \det}^{+} (\Phi),$$
où ${\det}^{+} (\Phi)$ est le produit des valeurs propres de module $>1$.
De plus, la mesure de Haar est l'unique mesure borélienne d'entropie maximale.

Considérons un élément $\phi=\sum_{h\in G} \phi_{h} h$ dans l'anneau entier $\mathbf{Z}[G]$ d'un groupe dénombrable $G$ et le quotient $\mathbf{Z}[G]/\mathbf{Z}[G] \phi$ par l'idéal à gauche engendré par $\phi$. Le dual de Pontryagin de ce groupe abélien discret est un groupe abélien compact muni d'une action continue par automorphismes de groupe, induite par la multiplication à gauche
\[G\curvearrowright X_{\phi}\overset{\mathrm{def}}{=}\widehat{(\mathbf{Z}[G]/\mathbf{Z}[G] \phi)}.\]
C'est le fermé $G$-invariant du décalage de Bernoulli $\widehat{\mathbf{Z}[G]}=\left( \mathbf{R}/\mathbf{Z}\right)^{G}$ formé des suites
\[X_{\phi}=\{(x_{g})_{g}\in \left( \mathbf{R}/\mathbf{Z}\right)^{G} \vert \sum_{h\in G} \phi_{h} x_{gh}=0, \forall g\in G\}.\]

On a toute une série de travaux qui permettent d'exprimer l'entropie topologique d'une telle action.
S.~Juzvinski{\u\i} \cite{Jusvinskii-1967-ent-gp-endom} dans le cas de $G=\mathbf{Z}$, l'exprime comme logarithme du produit des racines de module $\geqslant 1$ de  $\Phi$.

Le lemme de Mahler  \cite{Mahler-1960-app-Jensen-Form,Mahler-1962-ineq-sev-var} permet d'interpréter un tel produit en analyse complexe comme une intégrale (\og mesure de Mahler logarithmique \fg).  D.~Lind, K.~ Schmidt et T.~Ward \cite{Lind-Schmidt-Ward-1990-Mahler-mes} exprimeront l'entropie dans le cas de $G=\mathbf{Z}^p$ en ces termes, considérant $\Phi$ comme polynôme de Laurent
\begin{equation}
h_{\mathrm{top}}(\mathbf{Z}^p\curvearrowright X_{\Phi})=\int_{(\mathbf{R}/\mathbf{Z})^p} \log\vert \Phi(\mathrm{e}^{2\pi \mathrm{i}\theta})\vert \, d\theta.
\end{equation}
C.~Deninger \cite{Deninger-2006-Fuglede-Kadison-det} observe que ces quantités admettent des généralisations dans le cadre non commutatif à l'aide du déterminant de Fuglede-Kadison. Il s'agit d'un objet d'analyse fonctionnelle, concocté à l'aide du calcul fonctionnel et de la trace de von Neumann défini sur des opérateurs $u$ de 
$G$-modules de Hilbert et qui joue le rôle du déterminant positif classique ${\det}^{+}$~:
\[ {\det}^{+}_{\mathrm{vN}(G)} (u) \overset{\mathrm{def}}{=} \exp\left(\int_{]0,\infty[} \log(t) \, d\lambda_{\vert u\vert} (t)\right).\]
où $\lambda_{\vert u\vert}$ représente la fonction de densité spectrale de l'opérateur ${\vert u\vert}$.
Le domaine de validité de l'égalité entre l'entropie topologique et ce déterminant positif est peu à peu étendu, sous des hypothèses plus ou moins fortes sur $\phi$ (positivité, diverses formes d'inversibilité,...), à des classes de plus en plus grandes de groupes (croissance polynomiale, moyennable résiduellement fini, ...)
\cite{Deninger-2006-Fuglede-Kadison-det, Deninger-Schmidt-2007-exp-alg-act-entropy,Li-2012-Fuglede-Kadison-det, Li-Thom-2014-entropy-det-L2-torsion}
jusqu'à atteindre une forme optimale\footnote{En effet, puisque $\det^{+}_{\mathrm{vN}(G)} (\phi)$ est fini, une égalité ne sera envisageable que lorsque $h_{\mathrm{top}}^{\Sigma}(G\curvearrowright X_{\phi})<\infty$.}, due à B.~Hayes.
Ces travaux montrent également que la mesure de Haar maximise l'entropie mesurée.
\begin{theo}
 \cite[Th. 1.1]{Hayes-2014-Fuglede-Kadison-sofic-ent}
 Soit $G$ un groupe dénombrable sofique et $\Sigma$ une approximation sofique de $G$.
 Soit $\phi\in \mathrm{Mat}_{p,q} (\mathbf{Z}[G])$ et $G\curvearrowright X_{\phi}$ l'action algébrique associée.
 On a les propriétés suivantes.
 \begin{itemize}
\item [(i)] L'entropie sofique topologique $h_{\mathrm{top}}^{\Sigma}(G\curvearrowright X_{\phi})$ est finie si et seulement si $\phi$ est injective comme opérateur $\ell^{2}(G)^{\oplus q}\to \ell^{2}(G)^{\oplus p}$.
\end{itemize}
Supposons que $\phi$ est injective comme opérateur $\ell^{2}(G)^{\oplus q}\to \ell^{2}(G)^{\oplus p}$.
 \begin{itemize}
\item[(ii)] Si $p=q$,
alors 
$h_{\mathrm{top}}^{\Sigma}(G\curvearrowright X_{\phi})=h_{\mathrm{mes}}^{\Sigma}(G\curvearrowright X_{\phi}, \mathrm{Haar})=\log \det^{+}_{\mathrm{vN}(G)} (\phi)$.
\item[(iii)] Si $p\not =q$, alors 
$h_{\mathrm{top}}^{\Sigma}(G\curvearrowright X_{\phi})\leqslant h_{\mathrm{mes}}^{\Sigma}(G\curvearrowright X_{\phi}, \mathrm{Haar}) \leqslant \log \det^{+}_{\mathrm{vN}(G)} (\phi)$.
 \end{itemize}
\end{theo}

\subsection{Entropie d'Abért-Weiss}

Le calcul de l'entropie sofique mesurée des décalages de Bernoulli  (voir section~\ref{sect: ent sofique Bernoulli}) fait intervenir une mesure auxiliaire judicieusement choisie sur l'ensemble des partitions $K^{D_n}$ avant comptage, dans un esprit de \og grande déviation\fg.

Cette idée a conduit M.~Abért et B.~Weiss (communication personnelle de M.~Abért et annonce \cite{Weiss2015-survey-sofic-ent})
à une approche un peu différente de l'entropie sofique mesurée, dont on donne ici quelques éléments.

Soit $G\curvearrowright (X,\mu)$ une action p.m.p. du groupe sofique $G$ et $\alpha \colon X\to K$ une partition mesurée finie génératrice.

Pour chaque partie finie $F\subset G$, on cherche à imiter la mesure poussée en avant $\mu_{F}=\alpha^{\vee F}_{*} \mu$ sur $K^{F}$.
Ce qu'on a fait jusqu'ici consistait, partant d'une application 
\begin{equation}
\Upsilon \colon  \left(\begin{array}{ccl}
D\times K^{D} &\to &K^F 
\\
({v}, a) &\mapsto & (a(\sigma(f)\cdot {v}))_{f\in F}
\end{array}\right)
\end{equation}
à considérer, pour chaque partition $a\in K^{D}$ de $D$, la mesure poussée en avant 
$\Upsilon(\cdot, a)_{*}\mathbf{u} $
sur $K^{F}$ de la mesure uniforme sur $D$, puis à compter les bons $a$ (ceux pour lesquels cette mesure est proche de $\mu_{F}$).

La démarche de M.~Abért et B.~Weiss consiste à considérer des mesures de probabilité $\nu\in M(K^{D})$ sur l'ensemble fini $K^{D}$. Elles ont chacune une certaine entropie de Shannon $H(\nu)$ qui aura tendance à croître sous-linéairement en $\vert D\vert$, d'où la pertinence d'une normalisation $\frac{H(\nu)}{\vert D\vert}$. Pour chaque point $v\in D$, ils considèrent la mesure poussée en avant $\nu_{v,F}=\Upsilon(v, \cdot)_{*}\nu$. C'est une mesure sur $K^{D}$ qu'ils comparent avec $\mu_{F}$ pour la norme $\ell^1$, en moyenne sur $D$. Et ils définissent une notion d'entropie  
qui satisfait elle aussi les conditions (a) et (b). 
 \begin{equation}
h_{\mathrm{AW}}(G\curvearrowright X,\mu)\overset{\mathrm{def}}{=}\inf_{\epsilon > 0} \ \ \inf_{\substack{F \subset G\\F \text{ fini}}}
 \sup \frac{H(\nu)}{\vert D\vert},
 \tag{\textbf{Entropie selon Abért-Weiss}}
 \end{equation}
 où le \emph{supremum} est pris sur toutes les $(F,\epsilon)$-approximations sofiques $\Sigma=(\sigma \colon G\to \mathrm{Sym}(D))$ et pour toutes les mesures de probabilité $\nu\in M(K^{D})$ telles que $\frac{1}{\vert D\vert} \sum_{v\in D} \Vert \nu_{v,F}-\mu_{F}\Vert_1 <\epsilon$.
Il peuvent montrer que $h_{\mathrm{AW}}(G\curvearrowright X,\mu)\leqslant h_{\mathrm{mes}}^{\Sigma}(G\curvearrowright X,\mu)$ et L.~Bowen a exhibé des exemples\footnote{Qui encore une fois reposent sur la propriété $(\tau)$, voir section~\ref{sect: depend-approx-sofic}.} où on a une inégalité stricte.

\section{Remerciements}

Je veux exprimer ma gratitude à Alessandro Carderi et Mika\"el de la Salle dont l'aide m'a été très précieuse pour comprendre certaines références et pour la prépa\-ra\-tion et des relectures de ce texte.
Je tiens à remercier chaleureusement Mikl\'os Abért, Lewis Bowen, David Kerr et Brandon Seward qui m'ont initié à l'entropie sofique, m'ont fait part de leur vision, et ont répondu à mes nombreuses questions. 
David Kerr et Hanfeng Li  m'ont permis d'accéder à une version préliminaire du chapitre 9 \og {\em Entropy for actions of sofic groups} \fg\ de leur livre en préparation \cite{Kerr-Li-book}, je leur en suis reconnaissant. 
Merci également à toutes les personnes m’ayant signalé diverses coquilles, omissions ou imprécisions, notamment Pierre-Emmanuel Caprace, Véronique Chen Ai Ti, \'Etienne Ghys, Julien Melleray,
Romain Tessera ainsi qu'à toutes les autres personnes avec qui j'ai pu discuter de ces sujets, parmi lesquelles Tim Austin, Tullio Ceccherini-Silberstein, Grégory Miermont, Jean-François Quint,...

\thanks{Travail soutenu par le C.N.R.S. et  par le projet
  ANR-14-CE25-0004 GAMME.}

\bibliographystyle{alpha}

\begin{thebibliography}{}

\end{thebibliography}


\begin{thebibliography}{{Mey}15}

\bibitem[AE12]{Abert-Elek=profinite-act-2012}
M.~Ab{\'e}rt et G.~Elek.
\newblock Dynamical properties of profinite actions.
\newblock {\em Ergodic Theory Dynam. Systems}, 32 (6) : 1805--1835, 2012.

\bibitem[AKM65]{Adler-Konheim-McAndrew-1965-top-entropy}
R.~L. {Adler}, A.~G. {Konheim}, et M.~H. {McAndrew}.
\newblock {Topological entropy.}
\newblock {\em {Trans. Amer. Math. Soc.}}, 114 : 309--319, 1965.

\bibitem[AS16]{Seward-2016-Krieger-finite-th-Rokhlin-3}
A.~{Alpeev} et B.~{Seward}.
\newblock {Krieger's finite generator theorem for ergodic actions of countable
  groups III}.
\newblock {\em En préparation}, 2016.

\bibitem[{Aus}15]{Austin-2015-Add-prod-sofic}
T.~{Austin}.
\newblock {Additivity properties of sofic entropy and measures on model
  spaces}.
\newblock {\em ArXiv e-prints}, October 2015.

\bibitem[Bow71]{Bowen-1971-ent-top}
R.~Bowen.
\newblock Entropy for group endomorphisms and homogeneous spaces.
\newblock {\em Trans. Amer. Math. Soc.}, 153 : 401--414, 1971.

\bibitem[Bow10a]{Bowen-2010-f-invariant}
L.~Bowen.
\newblock A measure-conjugacy invariant for free group actions.
\newblock {\em Ann. of Math. (2)}, 171 (2) : 1387--1400, 2010.

\bibitem[Bow10b]{Bowen-2010-mes-conj-inv-sofic}
L.~Bowen.
\newblock Measure conjugacy invariants for actions of countable sofic groups.
\newblock {\em J. Amer. Math. Soc.}, 23 (1) : 217--245, 2010.

\bibitem[Bow11]{Bowen-2011=weak-isom-Bernoulli}
L.~Bowen.
\newblock Weak isomorphisms between {B}ernoulli shifts.
\newblock {\em Israel J. Math.}, 183 : 93--102, 2011.

\bibitem[Bow12a]{Bowen-2012-sofic-ent-amenab}
L.~Bowen.
\newblock Sofic entropy and amenable groups.
\newblock {\em Ergodic Theory Dynam. Systems}, 32 (2) : 427--466, 2012.

\bibitem[Bow12b]{Bowen-2012-almost-Ornstein}
L.~Bowen.
\newblock Every countably infinite group is almost {O}rnstein.
\newblock In {\em Dynamical systems and group actions}, volume 567 of {\em
  Contemp. Math.}, pages 67--78. Amer. Math. Soc., Providence, RI, 2012.

\bibitem[{Car}15]{Carderi-2015-arxiv}
A.~{Carderi}.
\newblock {Ultraproducts, weak equivalence and sofic entropy}.
\newblock {\em ArXiv e-prints}, September 2015.

\bibitem[CL15]{Capraro-Lupini-LNM-sofic-hyperlin-gps}
V.~Capraro et M.~Lupini.
\newblock {\em Introduction to {S}ofic and hyperlinear groups and {C}onnes'
  embedding conjecture}, volume 2136 of {\em Lecture Notes in Mathematics}.
\newblock Springer, Cham, 2015.
\newblock With an appendix by Vladimir Pestov.

\bibitem[Con73]{Conze-1972}
J.-P. Conze.
\newblock Entropie d'un groupe ab\'elien de transformations.
\newblock {\em Z. Wahrscheinlichkeitstheorie und Verw. Gebiete}, 25 : 11--30,
  1972/73.

\bibitem[Den74]{Denker-1974}
M.~Denker.
\newblock Finite generators for ergodic, measure-preserving transformations.
\newblock {\em Z. Wahrscheinlichkeitstheorie und Verw. Gebiete}, 29 : 45--55,
  1974.

\bibitem[Den06]{Deninger-2006-Fuglede-Kadison-det}
C.~Deninger.
\newblock Fuglede-{K}adison determinants and entropy for actions of discrete
  amenable groups.
\newblock {\em J. Amer. Math. Soc.}, 19 (3) : 737--758 (electronic), 2006.

\bibitem[Din70]{Dinaburg-1970-announct-of-Dinaburg-1971}
E.~I. Dinaburg.
\newblock A correlation between topological entropy and metric entropy.
\newblock {\em Dokl. Akad. Nauk SSSR}, 190 : 19--22, 1970.

\bibitem[Din71]{Dinaburg-1971-conn-var-ent-charac}
E.~I. Dinaburg.
\newblock A connection between various entropy characterizations of dynamical
  systems.
\newblock {\em Izv. Akad. Nauk SSSR Ser. Mat.}, 35 : 324--366, 1971.

\bibitem[DS07]{Deninger-Schmidt-2007-exp-alg-act-entropy}
C.~Deninger et K.~Schmidt.
\newblock Expansive algebraic actions of discrete residually finite amenable
  groups and their entropy.
\newblock {\em Ergodic Theory Dynam. Systems}, 27 (3) : 769--786, 2007.

\bibitem[ES05]{2005=Elek-Szabo=hyperlinearity}
G.~Elek et E.~Szab{\'o}.
\newblock Hyperlinearity, essentially free actions and {$L^2$}-invariants.
  {T}he sofic property.
\newblock {\em Math. Ann.}, 332 (2) : 421--441, 2005.

\bibitem[ES11]{Elek-Szabo-2011-sofic-amalg-over-amenable}
G.~Elek et E.~Szab{\'o}.
\newblock Sofic representations of amenable groups.
\newblock {\em Proc. Amer. Math. Soc.}, 139 (12) : 4285--4291, 2011.

\bibitem[Gab00]{Gab00a}
{}D. Gaboriau.
\newblock Co\^ut des relations d'\'equivalence et des groupes.
\newblock {\em Invent. Math.}, 139 (1) : 41--98, 2000.

\bibitem[GK76]{Grillenberger-Krengel-1976-Krieger-gen-th}
C.~Grillenberger et U.~Krengel.
\newblock On marginal distributions and isomorphisms of stationary processes.
\newblock {\em Math. Z.}, 149 (2) : 131--154, 1976.

\bibitem[Goo69]{Goodwyn-1969-maj-ent-top-by-mes}
L.~W.~Goodwyn.
\newblock Topological entropy bounds measure-theoretic entropy.
\newblock {\em Proc. Amer. Math. Soc.}, 23 : 679--688, 1969.

\bibitem[Goo71]{Goodman-1971-princ-variat}
T.~N.~T. Goodman.
\newblock Relating topological entropy and measure entropy.
\newblock {\em Bull. London Math. Soc.}, 3 : 176--180, 1971.

\bibitem[Gro99]{Gromov-1999-symbolic-algebraic-varieties}
M.~Gromov.
\newblock Endomorphisms of symbolic algebraic varieties.
\newblock {\em J. Eur. Math. Soc. (JEMS)}, 1 (2) : 109--197, 1999.

\bibitem[GS15]{Gab-Seward-2015-arxiv}
D.~{Gaboriau} et B.~{Seward}.
\newblock {Cost, $\ell^2$-Betti numbers and the sofic entropy of some algebraic
  actions}.
\newblock {\em ArXiv e-prints}, September 2015.

\bibitem[{Hay}14]{Hayes-2014-Fuglede-Kadison-sofic-ent}
B.~{Hayes}.
\newblock {Fuglede-Kadison determinants and sofic entropy}.
\newblock {\em ArXiv e-prints}, February 2014.

\bibitem[Juz67]{Jusvinskii-1967-ent-gp-endom}
S.~A. Juzvinski{\u\i}.
\newblock Calculation of the entropy of a group-endomorphism.
\newblock {\em Sibirsk. Mat. \u Z.}, 8 : 230--239, 1967.

\bibitem[Kat07]{Katok-2007-50-yrs-entropy}
A.~Katok.
\newblock Fifty years of entropy in dynamics: 1958--2007.
\newblock {\em J. Mod. Dyn.}, 1 (4) : 545--596, 2007.

\bibitem[Ker13]{Kerr-2013-Sofic-meas-ent-via-finite-partitions}
D.~Kerr.
\newblock Sofic measure entropy via finite partitions.
\newblock {\em Groups Geom. Dyn.}, 7 (3) : 617--632, 2013.

\bibitem[KL11a]{Kerr-Li-2011-Variationnal-principle}
D.~Kerr et H.~Li.
\newblock Entropy and the variational principle for actions of sofic groups.
\newblock {\em Invent. Math.}, 186 (3) : 501--558, 2011.

\bibitem[KL11b]{Kerr-Li-2011-Bernoulli-infinite-entropy}
D.~Kerr et H.~Li.
\newblock Bernoulli actions and infinite entropy.
\newblock {\em Groups Geom. Dyn.}, 5 (3) : 663--672, 2011.

\bibitem[KL13a]{Kerr-Li-2013-sofic-amenabl-dyn-entrop}
D.~Kerr et H.~Li.
\newblock Soficity, amenability, and dynamical entropy.
\newblock {\em Amer. J. Math.}, 135 (3) : 721--761, 2013.

\bibitem[KL13b]{Kerr-Li-2013-ent-top-version-combinat}
D.~Kerr et H.~Li.
\newblock Combinatorial independence and sofic entropy.
\newblock {\em Commun. Math. Stat.}, 1 (2) : 213--257, 2013.

\bibitem[KL15]{Kerr-Li-book}
D.~Kerr et H.~Li.
\newblock {\em Ergodic Theory Independence and Dichotomies}.
\newblock In preparation, 2015.

\bibitem[Kol58]{Kolmogorov-1958-entropy}
A.~N. Kolmogorov.
\newblock A new metric invariant of transient dynamical systems and
  automorphisms in {L}ebesgue spaces.
\newblock {\em Dokl. Akad. Nauk SSSR (N.S.)}, 119 : 861--864, 1958.

\bibitem[Kol59]{Kolmogorov-1959-entropy-erratum}
A.~N. Kolmogorov.
\newblock Entropy per unit time as a metric invariant of automorphisms.
\newblock {\em Dokl. Akad. Nauk SSSR}, 124 : 754--755, 1959.

\bibitem[Kri70]{Krieger-1970-finite-generator}
W.~Krieger.
\newblock On entropy and generators of measure-preserving transformations.
\newblock {\em Trans. Amer. Math. Soc.}, 149 : 453--464, 1970.

\bibitem[Li12]{Li-2012-Fuglede-Kadison-det}
H.~Li.
\newblock Compact group automorphisms, addition formulas and
  {F}uglede-{K}adison determinants.
\newblock {\em Ann. of Math.  (2)}, 176 (1) : 303--347, 2012.

\bibitem[LSW90]{Lind-Schmidt-Ward-1990-Mahler-mes}
D.~Lind, K.~Schmidt, et T.~Ward.
\newblock Mahler measure and entropy for commuting automorphisms of compact
  groups.
\newblock {\em Invent. Math.}, 101 (3) : 593--629, 1990.

\bibitem[LT14]{Li-Thom-2014-entropy-det-L2-torsion}
H.~Li et A.~Thom.
\newblock Entropy, determinants, and {$L^2$}-torsion.
\newblock {\em J. Amer. Math. Soc.}, 27 (1) : 239--292, 2014.

\bibitem[Lub94]{Lubotzky=book-exp-graph-1994}
A.~Lubotzky.
\newblock {\em Discrete groups, expanding graphs and invariant measures},
  volume 125 of {\em Progress in Mathematics}.
\newblock Birkh\"auser Verlag, Basel, 1994.
\newblock With an appendix by Jonathan D. Rogawski.

\bibitem[L{\"u}c94]{Luc94b}
{}W. L{\"u}ck.
\newblock Approximating ${L}\sp 2$-invariants by their finite-dimensional
  analogues.
\newblock {\em Geom. Funct. Anal.}, 4 (4) : 455--481, 1994.

\bibitem[Mah60]{Mahler-1960-app-Jensen-Form}
K.~Mahler.
\newblock An application of {J}ensen's formula to polynomials.
\newblock {\em Mathematika}, 7 : 98--100, 1960.

\bibitem[Mah62]{Mahler-1962-ineq-sev-var}
K.~Mahler.
\newblock On some inequalities for polynomials in several variables.
\newblock {\em J. London Math. Soc.}, 37 : 341--344, 1962.

\bibitem[{Mey}15]{Meyerovitch-2015-ent-stabilizer}
T.~{Meyerovitch}.
\newblock {Positive sofic entropy implies finite stabilizer}.
\newblock {\em ArXiv e-prints}, April 2015.

\bibitem[Orn70]{Ornstein-1970-entrop-shift-isom}
D.~Ornstein.
\newblock Bernoulli shifts with the same entropy are isomorphic.
\newblock {\em Advances in Math.}, 4 : 337--352  (1970), 1970.

\bibitem[Orn13]{Ornstein-2013-survey-Bernoulli}
D.~Ornstein.
\newblock Newton's laws and coin tossing.
\newblock {\em Notices Amer. Math. Soc.}, 60 (4) : 450--459, 2013.

\bibitem[OW87]{OW87}
{}D. Ornstein et {}B. Weiss.
\newblock Entropy and isomorphism theorems for actions of amenable groups.
\newblock {\em J. Analyse Math.}, 48 : 1--141, 1987.

\bibitem[P{\u{a}}u11]{Paunescu-2011-sofic-act-equiv-rel}
L.~P{\u{a}}unescu.
\newblock On sofic actions and equivalence relations.
\newblock {\em J. Funct. Anal.}, 261 (9) : 2461--2485, 2011.

\bibitem[Pes08]{Pestov-2008-sofic-gps-survey}
V.~G. Pestov.
\newblock Hyperlinear and sofic groups:  a brief guide.
\newblock {\em Bull. Symbolic Logic}, 14 (4) : 449--480, 2008.

\bibitem[Pop06]{Popa-1-cohomology-2006}
S.~Popa.
\newblock Some computations of 1-cohomology groups and construction of
  non-orbit-equivalent actions.
\newblock {\em J. Inst. Math. Jussieu}, 5 (2) : 309--332, 2006.

\bibitem[PS07]{PS07}
S.~Popa et R.~Sasyk.
\newblock On the cohomology of {B}ernoulli actions.
\newblock {\em Ergodic Theory Dynam. Systems}, 27 (1) : 241--251, 2007.

\bibitem[Roh63]{Rohlin-1963-generators}
V.~A. Rohlin.
\newblock Generators in ergodic theory.
\newblock {\em Vestnik Leningrad. Univ.}, 18 (1) : 26--32, 1963.

\bibitem[Roh67]{Rohlin-1967-lectures-entropy}
V.~A. Rohlin.
\newblock Lectures on the entropy theory of transformations with invariant
  measure.
\newblock {\em Uspehi Mat. Nauk}, 22 (5 (137)) : 3--56, 1967.


\bibitem[San57]{Sanov-1957-large-dev}
I.~N. Sanov.
\newblock On the probability of large deviations of random magnitudes.
\newblock {\em Mat. Sb. N. S.}, 42  (84) : 11--44, 1957.

\bibitem[{Sew}14]{Seward-2014-Krieger-finite-th-Rokhlin-1}
B.~{Seward}.
\newblock {Krieger's finite generator theorem for ergodic actions of countable
  groups I}.
\newblock {\em ArXiv e-prints}, May 2014.

\bibitem[{Sew}15]{Seward-2015-Krieger-finite-th-Rokhlin-2}
B.~{Seward}.
\newblock {Krieger's finite generator theorem for ergodic actions of countable
  groups II}.
\newblock {\em ArXiv e-prints}, January 2015.

\bibitem[Sha48]{Shannon-1948-entropy}
C.~E. Shannon.
\newblock A mathematical theory of communication.
\newblock {\em Bell System Tech. J.}, 27 : 379--423, 623--656, 1948.

\bibitem[Sin59]{Sinai-1959-entropy}
Y.~Sina{\u\i}.
\newblock On the concept of entropy for a dynamic system.
\newblock {\em Dokl. Akad. Nauk SSSR}, 124 : 768--771, 1959.

\bibitem[Sin62]{Sinai-1962-weak-isom}
Ja.~G. Sina{\u\i}.
\newblock A weak isomorphism of transformations with invariant measure.
\newblock {\em Dokl. Akad. Nauk SSSR}, 147 : 797--800, 1962.

\bibitem[ST14]{Seward-Tucker-Drob-2014-arxiv}
B.~{Seward} et R.~D. {Tucker-Drob}.
\newblock {Borel structurability on the 2-shift of a countable group}.
\newblock {\em ArXiv e-prints}, February 2014.

\bibitem[Ste75]{Stepin-1975}
A.~M. Stepin.
\newblock Bernoulli shifts on groups.
\newblock {\em Dokl. Akad. Nauk SSSR}, 223 (2) : 300--302, 1975.

\bibitem[Tho08]{Thom=Diophantine-approx=2008}
A.~Thom.
\newblock Sofic groups and {D}iophantine approximation.
\newblock {\em Comm. Pure Appl. Math.}, 61 (8) : 1155--1171, 2008.

\bibitem[Wei00]{Weiss-2000-sofic-gp}
B.~Weiss.
\newblock Sofic groups and dynamical systems.
\newblock {\em Sankhy\=a Ser. A}, 62 (3) : 350--359, 2000.
\newblock Ergodic theory and harmonic analysis  (Mumbai, 1999).

\bibitem[{Wei}15]{Weiss2015-survey-sofic-ent}
B.~{Weiss}.
\newblock {Entropy and actions of sofic groups}.
\newblock {\em {Discrete and Continuous Dynamical Systems - Series B}},
  20 (10) : 3375--3383, 2015.

\end{thebibliography}

\bigskip
\nobreak
\noindent \textsc{Damien Gaboriau
\\
C.N.R.S., 
\\
\'Ecole normale supérieure de Lyon,
\\
UMPA, UMR  5669,
\\
69364 Lyon
cedex 7, FRANCE}
\\
\noindent \texttt{gaboriau@ens-lyon.fr}
\end{document}